\documentclass[11pt]{amsart}
\usepackage{graphicx}

\usepackage{lmodern}			
\usepackage[T1]{fontenc}		
\usepackage[utf8]{inputenc}		
\usepackage{lastpage}			
\usepackage{indentfirst}		
\usepackage{color,xcolor}		
\usepackage{graphicx}			
\usepackage{microtype} 			
\usepackage[displaymath, mathlines]{lineno}

\usepackage[subentrycounter,seeautonumberlist,nonumberlist=true]{glossaries}

\usepackage{lipsum}			
\usepackage[portuguese,onelanguage]{algorithm2e}	
 \usepackage{amsbsy}			
 \usepackage{amscd}			
 \usepackage{amsfonts}			
 \usepackage{amsmath}			
 \usepackage{amssymb}			
 \usepackage{amstext}			
\usepackage{amsthm}			
\usepackage{hyperref}			
\usepackage{cleveref}			
\usepackage{dsfont}			
 \usepackage{ifthen}			
\usepackage{listings}           	
 \usepackage{lscape}             	
 \usepackage{mathrsfs}			
 \usepackage{pdfpages}           	

 \usepackage{verbatim}

\usepackage{latexsym}
\usepackage{amsthm}
\usepackage{times}
\usepackage{graphicx}
\usepackage{graphics}
\usepackage{hyperref}         
\usepackage{fancybox}
\usepackage{fancyhdr}
\usepackage{ae}
\usepackage[all]{xy}

\newcommand{\R}{\mathbb{R}}

\newcommand{\sen}{\text{sen}}

\newcommand{\lie}{\mathfrak{g}}
\newcommand{\Ad}{\text{Ad}}
\newcommand{\cqd}{\begin{flushright}$\Box$\end{flushright}}
\newcommand{\m}{\mathfrak{m}}
\newcommand{\sub}{\mathfrak{k}}
\newcommand{\h}{\mathfrak{h}}
\newcommand{\p}{\mathfrak{p}}
\newcommand{\q}{\mathfrak{q}}
\newcommand{\cvg}{\text{\bf g}_t}
\newcommand{\cvh}{\text{\bf h}_t}
\newcommand{\cvk}{\text{\bf k}_t}
\newcommand{\cvm}{\text{\bf m}_t}
\newcommand{\cvn}{\text{\bf n}_t}

\numberwithin{equation}{section}
\newtheorem{theorem}{Theorem}[section]
\newtheorem{proposition}[theorem]{Proposition}
\newtheorem{lemma}[theorem]{Lemma}
\newtheorem{corollary}[theorem]{Corollary}
\theoremstyle{definition}
\newtheorem{definition}[theorem]{Definition}

\theoremstyle{remark}
\newtheorem{remark}[theorem]{Remark}

\definecolor{verde}{rgb}{0,0.5,0}

\begin{document}

\title[Yamabe Problem on Aloff-Wallach Spaces]{Bifurcation and Local Rigidity of Homogeneous Solutions to the Yamabe Problem on Aloff-Wallach Spaces}

\author{Lino Grama}
\address{Instituto de Matemática, Estatística e Computação Científica, Universidade Estadual de Campinas,
Campinas, SP}
\email{linograma@gmail.com}
\author{Kennerson N. S. Lima}
\address{Unidade Acadêmica de Matemática, Universidade Federal de Campina Grande, 
Paraíba, PB}
\email{kennerson@mat.ufcg.edu.br}

\begin{abstract}
We construct 1-parameter families of well-known solutions to the Yamabe problem defined on Aloff-Wallach Spaces to determine bifurcation instants for these homogeneous spaces by examining changes in the Morse index of these metrics as the parameter varies over $]0,\infty[$. A bifurcation point for such families is an accumulation point of other solutions to the Yamabe problem, while a local rigidity point is an isolated solution of this problem, i.e., it is not a bifurcation point.

\end{abstract}


\maketitle

\section{Introduction}
Given a compact orientable Riemannian manifold $(M,g)$ with dimension $m\geq3$, the {\it Yamabe problem} concerns the existence of constant scalar curvature metrics on $M$ conformal to $g$. Solutions to this problem, called {\it Yamabe metrics}, can be characterized variationally as critical points of the Hilbert-Einstein functional 
\begin{equation*}
\mathcal{A}(g)=\dfrac{1}{\textrm{Vol}(g)^\frac{m-2}{m}}\int_M\textrm{scal}(g)\textrm{vol}_g, 
\end{equation*}
restricted to the set $[g]$ of metrics conformal to $g$. The existence of such solutions is consequence of the successive works of Yamabe \cite{yamabe}, Trudinger \cite{trundinger}, Aubin \cite{aubin} and Schoen \cite{schoen}.


In this paper, we focus on studying the Yamabe problem in the context of homogeneous manifolds $G/H$ and homogeneous $G$-invariant metrics. The investigation of the invariant version of the Yamabe problem has gained significant attention, with notable progress being made recently; see for instance,  \cite{pacificjournal}, \cite{bettiol} \cite{knsl}.

One of the features of the invariant version of the Yamabe problem is that the families of metrics we study are formed by homogeneous metrics, which are trivial solutions to the Yamabe problem. In this context, natural questions about \textit{bifurcation theory} arise.

A bifurcation point of a family of solutions to the Yamabe problem is an accumulation point of other solutions to this problem conformal to elements of the family. 

The bifurcation theory applied here consists of finding bifurcation instants for a 1-parameter family of homogeneous metrics $g_t, t>0$. For this, we  analyze the occurrence of {\it jump} in the Morse index of $g_t$ by using strongly expressions of the eigenvalues of the Laplacian defined on $(M,g_t)$ and formulae for the scalar curvatures of the canonical variations $g_t$.  An interesting aspect of the bifurcation phenomena we are studying is occurrence of {\it break of symmetry} at all bifurcation points, which in our situation is characterized by the non-homogeneity of these new Yamabe metrics, as well as in the cases of \cite{pacificjournal}, \cite{bettiol} and \cite{knsl}. In the case of product metrics discussed in the work \cite{piccioneproduct}, the break of symmetry is characterized by the fact that the bifurcation points are not product metrics.

Before we state our results, let us establish the setup in which we will work in this paper. From a {\it Riemannian submersion with totally geodesic fibers} it can be construct a 1-parameter family of other Riemannian submersions of the same type by scaling the original metric of the total space in the direction of the fibers. This family is called {\it canonical variation} and a generalization of this is in particular the {\it Cheeger deformations}. Special cases of Riemannian submersion with totally geodesic fibers are the so-called {\it homogeneous fibrations}.


The \textit{Aloff-Wallach spaces} are a family of homogeneous spaces described as the total space of the homogeneous fibration 
\begin{center}
$U(2)/S^1_{k,l} \, \cdots \, W^7_{k,l}=SU(3)/S^1_{k,l} \rightarrow SU(3)/U(2)=\mathbb{CP}^{2},$
\end{center}  
where $S^1_{k,l}=\textrm{diag}(e^{2\pi ik\theta},e^{2\pi il\theta},e^{2\pi im\theta})\subset U(2)$, $\mathbb{Z}\ni k,l\geq 0$, $k+l+m=0$, and $\gcd(k,l)=1$.

This family of spaces was introduced in \cite{AW} and constitutes an important family of homogeneous 7-manifolds, serving as the groundwork for several phenomena in geometry. Notable examples include metrics of positive sectional curvature \cite{berger}, classification of invariant Einstein metrics \cite{WZ}, and recent studies in invariant $G_2$ structures and gauge theory \cite{G2}, \cite{Goncalo}.

We equip each total space above with a {\it normal homogeneous metric} and deform these metrics by shrinking the fibers by a factor $t^2$, $t>0$ getting the {\em canonical variations} $(W^7_{k,l},\cvg)$. Recall that a normal homogeneous metric on $W^7_{k,l}=SU(3)/S^1_{k,l}$ is obtained from the restriction to the tangent space (isotropy representation) of a bi-invariant inner product on the Lie algebra $\mathfrak{su}(3)$ . 


A {\it degeneracy point} for a given canonical variation $g_t, t>0$, is a degenerate critical point $g_{t_{\ast}}$ for the Hilbert-Einstein functional, at some $t_{\ast}>0$. For this canonical variation every bifurcation point is a degeneracy point, but in general not all degeneracy points are bifurcation points. Theorem 4.1 in \cite{otoba} can be applied to our canonical variations to study if the degeneracy points are bifurcation points.  Alternatively, we showed that the Morse index of $\cvg$, changes as $t$ crosses {\it degeneracy instants} $t_{\ast}>0$ using the notion of {\it compensation of eigenvalues (or multiplicities)} to prove existence of new solutions to the Yamabe problem accumulating at $\text{\bf g}_{t_{\ast}}$. Such instants $t_{\ast}$ are called {\it bifurcation instants} and $\text{\bf g}_{t_{\ast}}$ is a bifurcation point for the canonical variations $\cvg$.

Regarding the study of bifurcation in this set up, Bettiol and Piccione proved in \cite{pacificjournal} that for every homogeneous fibration there exists a sequence $\{t_q\}_{q\in\mathbb{N}}$ in the interval $]0,1[$, converging to $0$, such that each $g_{t_q}$ is a bifurcation point for the corresponding canonical variation $\cvg$. Recently, Otoba and Petean in \cite{otoba} studied bifurcation for the constant scalar curvature equation along a 1-parameter family of Riemannian metrics on the total space of a harmonic Riemannian submersion and provide an existence theorem for bifurcation points. From the work of Otoba and Petean follows that the existence of infinitely many bifurcation points and the fact that the sets of such points are discrete, for each canonical variation of any Riemannian submersion with totally geodesic fibers. However, it is important to point out that there are no general methods that would allow us to exhibit explicitly all bifurcation points of a given canonical variation neither its Morse index. This was done only in a few particular cases, actually for {\it Berger spheres} by Bettiol and Piccione in \cite{bettiol} and for {\em full flag manifolds} by  the authors in the work \cite{knsl}. 

In this present work, our main contribution is to give the precise location of all bifurcation points, namely for the canonical variations above mentioned, in a very explicit way, applying Lie theoretic methods { \em for the Yamabe problem in Aloff-Wallach spaces}. This approach also enabled us to calculate the first positive eigenvalue of these canonical variations as well as their respective Morse index. From Proposition \ref{criteriorig} we also determine the local rigidity instants for $\cvg$.

It is worth pointing out that the analysis of the invariant Yamabe problem in Aloff-Wallach spaces is more involved compared to, for instance, flag manifolds, due to the rich geometry of such spaces. One remarkable feature of this family of spaces is that the isotropy representation has equivalent components, and therefore the description of invariant tensors needs to be carried out more carefully.

Our results also can be understood from the view point of dynamical systems, where {\it bifurcation} means a topological or qualitative change in the structure of the set of fixed points of a 1-parameter family of systems when we vary this parameter. Critical points of the Hilbert-Einstein functional in a conformal class $[g]$ are fixed points of the so-called {\it Yamabe flow}, the corresponding $L^2$-gradient flow of the Hilbert-Einstein functional, which gives a dynamical system in this conformal class. Hence, the bifurcation results above mentioned can be interpreted as a local change in the set of fixed points of the Yamabe flow near homogeneous metrics (which are always fixed points) when varying the conformal class $[g]$ with $g$ in one of the families $\cvg$. Then an interesting question would be to study the dynamics near these new fixed points.


In the bifurcation results mentioned, we analyze the second variation of the Hilbert-Einstein functional at the homogeneous metric $\cvg$, defined on the total space $(W^7_{k,l},\cvg)$. Given the second variation formula for this functional, such analysis is equivalent to comparing the eigenvalues $\lambda^{kj}(t)$ (see expression \eqref{kj}, Corollary \ref{deltat}) of the Laplacian $\Delta_t$ of $\cvg$ with the scalar curvatures $\text{scal}(\cvg)$. More precisely, a critical point $\cvg$ is degenerate if and only if $\text{scal}(t)\neq 0$ and $\dfrac{\text{scal}(t)}{m-1}$ is an eigenvalue of $\Delta_t.$ 

The spectrum of the Laplacian $\Delta_t$ of the canonical variations of Riemannian submersions with totally geodesic fibers is well-understood. Roughly, it consists of linear combinations (that depends on $t$) of eigenvalues of the original metric on the total space with eigenvalues of the fibers. In particular, we compute formulae for the first positive eigenvalues of $\Delta_{\cvg}$. Combining this knowledge of the spectra of $\Delta_{\cvg}$ with the formula for the scalar curvature $\text{scal}(\cvg)$ we are able to identify all degeneracy instants for $\cvg$ and prove existence of bifurcation at all degeneracy instants in the interval $]0,\infty[$. The local rigidity instants for these canonical variations are also determined for all other $t>0$. We also compute the Morse index of $\cvg$, for  all $t>0$.

In the study of the spectra both of Aloff-Wallach spaces endowed with a normal homogeneous metric and the spectrum of the Laplacian on compact isotropy irreducible Hermitian symmetric space $G/H$ we used the description presented in \cite{Urakawa} and \cite{sugiura}. It is important to point out that the homogeneous metric on $\cvg$ is not a normal homogeneous metric on the respective Aloff-Wallach space for $0t>0$ and the spectra of its Laplacian $\Delta_t$ is indeterminate in general.  

Using our results on the existence of infinitely many bifurcations points for the families $\cvg$, we obtain in Section \ref{Section3.3} for this family existence of a subset $\mathcal{G}\subset ]0,\infty[$, accumulating at $0$ and such that for each $t\in\mathcal{G}$, there are at least $3$ solutions to the Yamabe problem in each conformal class $[\cvg]$. 

Our main results are summarized as follows (see Section 5 for further details).

{\bf Theorem A} (=Theorem \ref{teo-princ01}). {\it Let $\cvg$ be the canonical variation on $W^7_{k,l}$, with $(k,l)\neq (1,0), (1,1)$, and take $$b=\sqrt{\dfrac{-96\gamma+\sqrt{(96\gamma)^2+4(16\gamma-3k^2)(10\gamma-3k^2)}}{2(10\gamma-3k^2)}}.$$ Thus, the degeneracy instants for $\cvg$ in $]0,1[$ form a decreasing sequence $\{t^{\text{\bf g}}_{q}\}\subset \, ]0,b]$ such that $t^{\text{\bf g}}_{q}\rightarrow 0$ as $q\rightarrow 0$, with $t^{\text{\bf g}}_1=b$ and for $q>1$, 
\begin{equation}
t^{\text{\bf g}}_{q}=\dfrac{1}{\sqrt{2(10\gamma-3k^2})}\sqrt{48\gamma(1-2q -q^2)+\sqrt{(48\gamma)^2(1-2q -q^2)^2-4(3k^2-10\gamma)(16\gamma -3k^2)}}, \label{seqsu}
\end{equation}
$0<l<k$, $k,l$ relatively prime and $\gamma=k^2+kl+l^2.$ }

The similar results for the manifolds $W^7_{1,0}$, $W^7_{1,1}$ are established as follows.

{\bf Theorem B} (=Theorem \ref{teo-princ-02}). {\it Let $\cvh$ be the canonical variation on $W^7_{1,0}$ and take $b=2.$ Thus, the degeneracy instants for $\cvh$ in $]0,b]$ form a decreasing sequence $\{t^{\text{\bf h}}_{q}\}\subset \, ]0,b]$ such that $t^{\text{\bf h}}_{q}\rightarrow 0$ as $q\rightarrow 0$, with $t^{\text{\bf h}}_1=b$ and for $q>1$, 
\begin{equation}
t^{\text{\bf h}}_{q}=\frac{\sqrt{24-4q-2q^2+\sqrt{\left(-2q^2-4q+24\right)^2+160}}}{\sqrt{10}}. \label{seqso1}
\end{equation}
}

{\bf Theorem C} (=Theorem \ref{teo-princ03}). {\it Let $\cvk$ be the canonical variation on $W^7_{1,1}$ and take $b=\sqrt{3+\sqrt{10}}.$ Thus, the degeneracy instants for $\cvk$ in $]0,b]$ form a decreasing sequence $\{t^{\text{\bf k}}_{q}\}\subset \, ]0,b]$ such that $t^{\text{\bf k}}_{q}\rightarrow 0$ as $q\rightarrow 0$, with $t^{\text{\bf k}}_1=b$ and for $q>1$, 
\begin{equation}
t^{\text{\bf k}}_{q}=\frac{\sqrt{24-4q-2q^2+\sqrt{\left(-2q^2-4q+24\right)^2+36}}}{\sqrt{6}}. \label{seqso1}
\end{equation}
}

This paper is organized as follows. In Section \ref{Section1} we introduce the variational characterization of the Yamabe problem and the basic notions of bifurcation and local rigidity. In Section \ref{AW} we recall the definitions of Aloff-Wallach spaces, invariant metrics and isotropy representation. The study of the Laplacian of the Riemannian submersions with totally geodesic fibers, the definition of the canonical variations and the description of their spectra, as well as the construction of the canonical variations $\cvg$, $\cvh$, $\cvk$ are studied in Section \ref{cvariation}. In Section \ref{section5} we obtain formulae for scalar curvature of $\cvg$, $\cvh$, $\cvk$. We also prove in this Section the main results of the work on bifurcation and local rigidity of solutions to the Yamabe problem, for the families $\cvg$, $\cvh$, $\cvk$. We also determine the Morse index of $\cvg$, $\cvh$, $\cvk$, for $t>0$. The last section of this chapter contains multiplicity of solutions to the Yamabe problem for all the families above mentioned.

\section{Variational Setup For the Yamabe Problem}\label{Section1}

Before we present the notions of bifurcation and local rigidity instants, we will describe the set ${\mathcal{R}}^k(M)$ of all {\it $C^k$ Riemannian metrics} on $M$, $m=\dim M\geq 3$, and some properties of the Hilbert-Einstein functional. For further references and details of the concepts and properties presented here see, e.g, \cite{bettiol}, \cite{piccioneproduct} and \cite{schoen}.
	
Let $g_R$ be a fixed auxiliary Riemannian metric on $M$; $g_R$ and its Levi-Civita connection $\nabla^R$ induce naturally Riemannian metrics and connections on all tensor bundles over $M$, respectively. For each $k\geq0$, denote by $\Gamma^k(TM^{\ast}\vee TM^{\ast})$ the space of symmetric $C^k$ sections of $TM^{\ast}\otimes TM^{\ast}$, i.e, symmetric $(0,2)$-tensors of class $C^k$ on $M$. This becomes a Banach space when equipped with the $C^k$ norm $$\left\|\tau\right\|_{C^k}=\displaystyle\max_{j=1,\ldots,k}\left(\displaystyle\max_{x\in M}\left\|(\nabla^R)^j\tau(x)\right\|_R\right),$$ where $\left\|\cdot\right\|_R$ denote the norm induced by $g_R$ on each appropriate space. 
	
Note that the set ${\mathcal{R}}^k(M)$ of all {\it $C^k$ Riemannian metrics} on $M$ is a open convex cone inside $(\Gamma^k(TM^{\ast}\vee TM^{\ast}),\left\|\cdot\right\|_{C^k})$ and, therefore, contractible. Thus, ${\mathcal{R}}^k(M)$ inherits a natural differential structure. We remark also that ${\mathcal{R}}^k(M)$ is an open set of a Banach space and, thus, we can identify its tangent space with the vector space $\Gamma^k(TM^{\ast}\vee TM^{\ast})$. From now on, we assume that $k\geq3.$
		
For each $g\in{\mathcal{R}}^k(M)$, denote by $\textrm{vol}_g$ the volume form on $M$ (we assume that $M$ is orientable); in this case, $L^2(M,\textrm{vol}_g)$ will denote the usual Hilbert space of real square integrable functions on $M$. Consider over ${\mathcal{R}}^k(M)$ the maps $$\textrm{scal}:{\mathcal{R}}^k(M)\longrightarrow C^{k-2}(M) \ \ \textrm{and} \ \ \textrm{Vol}:{\mathcal{R}}^k(M)\longrightarrow\R,$$ the {\it scalar curvature} and the {\it volume}, that for each Riemannian metric $g\in{\mathcal{R}}^k(M)$ associate, respectively, its scalar curvature $\textrm{scal}(g):M\longrightarrow\R$ and its volume $\textrm{Vol}(g)=\int_M{\textrm{vol}_g}$. Define the {\it Hilbert-Einstein functional} as the function $A:{\mathcal{R}}^k(M)\longrightarrow\R$ given by 
\begin{equation}
\mathcal{A}(g)=\dfrac{1}{\textrm{Vol}(g)^\frac{m-2}{m}}\int_M\textrm{scal}(g)\textrm{vol}_g. \label{HE}
\end{equation}
We will now establish the more appropriate regularity for our manifold of metrics and maps. Denote by ${\mathcal{R}}_{1}^k(M)=\textrm{Vol}^{-1}(1)$ the subset of ${\mathcal{R}}^k(M)$ consisting of metrics of volume one. Note that ${\mathcal{R}}_{1}^k(M)$ is a smooth embedded codimension 1 submanifold of ${\mathcal{R}}^k(M)$. 

For each $g\in{\mathcal{R}}^k(M)$, we denote by $[g]\subset {\mathcal{R}}^k(M)$ the set of metrics conformal to $g$. The space $\Gamma^k(TM^{\ast}\vee TM^{\ast})$ defined above induces a differential structure on each conformal class. 


Let $$[g]_{k,\alpha}=\left\{\phi g;\phi\in C^{k,\alpha}(M), \ \ \phi>0\right\}$$ be the {\it $C^{k,\alpha}(M)$ conformal class} of $g$, which can be identified with the open subset of $C^{k,\alpha}(M)$ formed by the positive functions, which allows $[g]_{k,\alpha}$ to have a natural differential structure. Working with this regularity, the necessary Fredholm's conditions for the second variation of the Hilbert-Einstein functional are satisfied. The set $${\mathcal{R}}_{1}^{k,\alpha}(M,g)={\mathcal{R}}_{1}^k(M)\cap[g]_{k,\alpha}$$is a smooth embedded codimension 1 submanifold  of $[g]_{k,\alpha}$. 

\begin{proposition} [\cite{piccioneproduct}] \label{hilbert}The restriction of the Hilbert-Einstein functional $\mathcal{A}$ to $\mathcal{R}_1^{k,\alpha}(M,g)$ is smooth, and its critical points are constant scalar curvature metrics in $\mathcal{R}_1^{k,\alpha}(M,g)$. 
\end{proposition}

According \cite{piccioneproduct} we also have that given a critical point $g_0\in \mathcal{R}_1^{k,\alpha}(M,g)$ of $\mathcal{A}$, the second variation of $\mathcal{A}$ at $g_0$ can be identified with the quadratic form $$\text{d}^2\mathcal{A}(g_0)(\psi,\psi)=\dfrac{m-2}{2}\int_M((m-1)\Delta_{g_0}\psi-\text{scal}(g_0)\psi)\psi\text{vol}_{g_0},$$ defined on the tangent space at $g_0=\phi g$, given by $$T_{g_0}{\mathcal{R}}_{1}^{k,\alpha}(M,g)\cong\left\{\psi\in C^{k,\alpha}(M);\int_M\frac{\psi}{\phi}\text{vol}_{g_0}=0\right\}.$$

\begin{remark} We are considering $\Delta_{g}=-\text{div}_g\circ\text{grad}_g$, the Laplacian operator of $(M,g)$, acting on $C^{\infty}(M)$. We observe that, given $\lambda\in\R^+$, one has $\Delta_{\lambda g}=\frac{1}{\lambda}\Delta_g$ and $\text{scal}(\lambda g)=\frac{1}{\lambda}\text{scal}(g)$. Therefore,  we can normalize metrics $g$ to have unit volume, without changing the spectral theory of the operator $\Delta_g-\dfrac{\text{scal}(g)}{m-1}\cdot\text{Id}$, $m=\dim M$.
\end{remark} 
 
Now, we will introduce the concepts of {\it bifurcation} and {\it local rigidity} for a 1-parameter family of solutions to the Yamabe problem. Let $$[a,b] \ni t \mapsto g_t\in {\mathcal{R}}^k(M), \ \ k\geq 3,$$ or simply $g_t$, denote a smooth 1-parameter family (smooth path) of Riemannian metrics on $M$, such that each $g_t$ has constant scalar curvature.
 
\begin{definition} An instant $t_{\ast}\in [a,b]$ is a \emph{bifurcation instant} for the family $g_t$ if there exists a sequence $\{t_q\}$ in $[a,b]$ that converges to $t_{\ast}$ and a sequence $\{g_q\}$ in ${\mathcal{R}}^k(M)$ of Riemannian metrics that converges to $g_{t_{\ast}}$ satisfying

\begin{itemize}
\item[(i)] $g_{t_q}\in [g_q]$, but $g_q\neq g_{t_q}$;
\item[(ii)] $\emph{Vol}(g_q)=\emph{Vol}(g_{t_q})$;
\item[(iii)]$\emph{scal}(g_q)$ is constant.
\end{itemize}
\end{definition}

If $t_{\ast}\in [a,b]$ is not a bifurcation instant, it is said that the family $g_t$ is {\it locally rigid} at $t_{\ast}$. More precisely, the family is locally rigid at $t_{\ast}\in [a,b]$ if there exists an open set $U\subset{\mathcal{R}}^k(M)$ containing $g_{t_{\ast}}$ such that if $g\in U$ is another metric with constant scalar curvature and there exists $t\in [a,b]$ with $g_t\in U$ and

\begin{itemize}
\item[\textrm{(a)}] $g\in [g_t]$;
\item[\textrm{(b)}] $\textrm{Vol}(g)=\textrm{Vol}(g_t)$,
\end{itemize} then $g=g_t$.


\begin{definition} An instant $t_{\ast}\in [a,b]$ is a \emph{degeneracy instant} for the family $g_t$ if $\emph{scal}(g_{t_{\ast}})\neq 0$ and $\frac{\emph{scal}(g_{t_{\ast}})}{m-1}$ is an eigenvalue of the Laplacian operator $\Delta_{t_{\ast}}$ of $g_{t_{\ast}}$. 
\end{definition}

\begin{remark} In face of the second variation expression of the Hilbert-Einstein functional the above definition is equivalent to the fact that $g_{t_{\ast}}$ being a degenerate critical point (in the usual sense of Morse theory) of the Hilbert-Einstein functional on ${\mathcal{R}}_{1}^{k,\alpha}(M,g)$.
\end{remark}

The {\it Morse index} $N(g_t)$ of $g_t$ is given by the number of positive eigenvalues of $\Delta_t$ counted with multiplicity less than $\frac{\text{scal}(g_{t_{\ast}})}{m-1}$. For each $t>0$, $N(g_t)$ is a non-negative integer number. In other words, $N(g_t)$ counts the number of directions which the functional decreases, since the second variation is negative definite.

We will now state a general criterion for classifying local rigidity instants for a 1-parameter family $g_t$ of solutions to the Yamabe problem.

\begin{proposition}[\emph{\cite{bettiol}}]\label{criteriorig} Let $g_t$ be a smooth path of metrics of class $C^k$, $k\geq 3$, such that $\emph{scal}(g_t)$ is constant for all $t\in [a,b]$, and let $\Delta_t$ be the Laplacian operator of $g_t$. If $t_{\ast}$ is not a degeneracy instant of $g_t$, then $g_t$ is locally rigid at $t_{\ast}$.
\end{proposition}

\begin{corollary}[\emph{\cite{bettiol}}] \label{criterioauto}Suppose that in addition to the hypotheses of Proposition \ref{criteriorig}, there exists an instant $t_{\ast}$ when $\frac{\emph{scal}(g_{t_{\ast}})}{m-1}$ is less than the first positive eigenvalue $\lambda_1(t_{\ast})$ of $\Delta_{t_{\ast}}$. Then $g_{t_{\ast}}$ is a local minimum for the Hilbert-Einstein functional in its conformal class. In particular, $g_t$ locally rigid at $t_{\ast}.$
\end{corollary}

\begin{remark}\emph{Note that, in fact, the Corollary \ref{criterioauto} is an immediate consequence from Proposition \ref{criteriorig}, since Morse index of the non-degenerate critical point $g_{t_{\ast}}$ is $N(g_{t_{\ast}})=0.$ In addition, we remark that, if $t_{\ast}$ is a bifurcation instant for $g_t$, then $t_{\ast}$ is necessarily a degeneracy instant for $g_t$. However, the reciprocal is not true in general.} 
\end{remark}

The next result provides a sufficient condition to determinate if a degeneracy instant is a bifurcation instant, given in terms of a jump in the Morse index when passing a degeneracy instant. 

\begin{proposition}\emph{\cite{bettiol}} \label{criterioindice}Let $g_t$ a smooth path of metrics of class $C^k$, $k\geq 3$, such that $\emph{scal}(g_t)$ is constant for all $t\in [a,b]$, $\Delta_t$ the Laplacian $g_t$ and $N(g_t)$ the Morse index of $g_t$. Assume that $a$ and $b$ are not degeneracy instants for $g_t$ and $N(g_a)\neq N(g_b)$. Then, there exists a bifurcation instant $t_{\ast}\in\, ]a,b[$ for the family $g_t$.
\end{proposition}

\section{Aloff-Wallach Spaces}\label{AW}
	
%
%
The Aloff-Wallach spaces are seven dimensional homogeneous spaces given by $W_{k,l}^7=SU(3)/SO(2)$ which can be described as $SU(3)/S^1_{k,l}$ where $S^1_{k,l}=\textrm{diag}(e^{2\pi ik\theta},e^{2\pi il\theta},e^{2\pi im\theta})\subset U(2),$ \linebreak $\mathbb{Z}\ni k,l\geq 0,\, k+l+m=0, \,\gcd(k,l)=1$, i.e, $SO(2)$ is seen as the circle embedding into $SU(3)$ by $$i_{l,k}:e^{2\pi il\theta}\mapsto \textrm{diag}(e^{2\pi ik\theta},e^{2\pi il\theta},e^{2\pi im\theta}). $$

Let us consider the following orthonormal basis of the Lie algebra $\mathfrak{su}(3)=\text{Lie}(SU(3))$, with respect to the inner product on $\mathfrak{su}(3)$ given by $Q(X,Y)=\langle X,Y \rangle =-\frac{1}{2}\text{ReTr}(XY)$:

 $$X_1 = (E_{12}-E_{21}), \quad X_2 = i(E_{12}-E_{21}),$$ 
 $$X_3 = (E_{13}-E_{31}), \quad X_4 = i(E_{13}-E_{31}),$$ 
 $$X_5 = (E_{23}-E_{32}), \quad X_6 = i(E_{23}-E_{32}),$$
 $$X_0=i\textrm{diag}[(k+2l),-(2k+l),(k-l)]/\sqrt{3\gamma},$$
 $$Z=i\textrm{diag}[(k,l,-(k+l)]/\sqrt{\gamma},$$
where $\gamma=k^2+kl+l^2$ and $E_{i,j}$ is a $3\times 3$ matrix with 1 on entry $(i,j)$ and 0 elsewhere. It follows that (see \cite{nikoronov}):

\begin{align*}
 [X_2,X_1]&=-(k+l)\sqrt{3}(\sqrt{\gamma})^{-1}X_0-(k-l)(\sqrt{\gamma})^{-1}Z  \\
 [X_3,X_1]&= X_5, & [X_4,X_1]&=X_6\\
 [X_5,X_1]&=X_3, & [X_6,X_1]&=X_4\\
 [X_0,X_1]&=3(k+l)(\sqrt{3\gamma})^{-1}X_2, & [Z,X_1]&=(k-l)(\sqrt{\gamma})^{-1}X_2\\
 [X_3,X_2]&=-X_6, & [X_4,X_2]&=X_5\\
 [X_5,X_2]&=-X_4, & [X_6,X_2]&=X_3\\
 [X_0,X_2]&=-(k+l)(\sqrt{3\gamma})^{-1}X_1, & [Z,X_2]&=-(k-l)(\sqrt{\gamma})^{-1}X_1 \\ 
 [X_4,X_3]&=-l\sqrt{3}(\sqrt{\gamma})^{-1}X_0-(2k+l)(\sqrt{\gamma})^{-1}Z  \\
 [X_5,X_3]&=X_1, & [X_6,X_3]&=-X_2\\
 [X_0,X_3]&=l\sqrt{3}(\sqrt{\gamma})^{-1}X_4, & [Z,X_3]&=(2k+l)(\sqrt{\gamma})^{-1}X_4  \\
 [X_5,X_4]&=X_2, & [X_6,X_4]&=X_1  \\
 [X_0,X_4]&=-l\sqrt{3}(\sqrt{\gamma})^{-1}X_3, & [Z,X_4]&=-(2k+l)(\sqrt{\gamma})^{-1}X_3  \\
 [X_6,X_5]&=k\sqrt{3}(\sqrt{\gamma})^{-1}X_0-(k+2l)(\sqrt{\gamma})^{-1}Z  \\
 [X_0,X_5]&=-k\sqrt{3}(\sqrt{\gamma})^{-1}X_6, & [Z,X_5]&=(k+2l)(\sqrt{\gamma})^{-1}X_6  \\
 [X_0,X_6]&=k\sqrt{3}(\sqrt{\gamma})^{-1}X_5, & [Z,X_6]&=-(k+2l)(\sqrt{\gamma})^{-1}X_5 \\
 [Z,X_0]&=0  \\
\end{align*}

\subsection{Isotropy Representation} \label{isotropia}  

From the notion of {\it isotropy representation} we can determine invariant metrics on certain homogeneous spaces as follows.
 
Let $G\times M\longrightarrow M$ be a differentiable and transitive action of a Lie group $G$ on the homogeneous space $(M,m)$ endowed with a $G$-invariant metric $m$. Given $x\in G$, let $K=G_x$ be the isotropy subgroup of $x$. The {\it isotropy representation} of $K$ is the homomorphism $g\in K\mapsto dg_x \in \textrm{Gl}(T_xM)$. Note that $m_{x_0}$ is a inner product on $T_{x_0}(G/K)$, invariant by such representation, where $x_0=1\cdot K$.

A homogeneous space $G/K$ is {\it reducible} if $G$ has a Lie algebra $\lie$ such that $$\lie=\mathfrak{k}\oplus\mathfrak{m},$$ with $\Ad(K)\mathfrak{m}\subset \mathfrak{m}$ and $\mathfrak{k}$ Lie algebra of $K$. If $K$ is compact, this decomposition always exists, namely, if we take $\mathfrak{m}=\mathfrak{k}^{\bot}$, $-B$-orthogonal complement to $\mathfrak{k}$ in $\lie$, where $B$ is the Cartan-Killing form of $\lie$.

If $G/K$ is reducible the isotropy representation of $K$ is equivalent to $\Ad |_K$, the restriction of the adjoint representation of $G$ to $K$, that is,  $$j(k)=\Ad(k)|_{\mathfrak{m}}, \forall k\in K.$$

A representation of a compact Lie group $K$ is always orthogonal (preserves inner product) on the representation space. We can conclude that every reductive homogeneous space $G/K$ has a $G$-invariant metric, since such a metric is completely determined by an inner product on the tangent space at the origin $T_{x_0}(G/K)$.  

We remark that the set of all $G$-invariant metrics on $G/K$ are in 1-1 correspondence the set of inner products $\left\langle, \right\rangle$ on $\mathfrak{m}$, invariant by $\Ad(k)$ on $\mathfrak{m}$, for each $k\in K$, that is, $$\left\langle \Ad(k)X, \Ad(k)Y \right\rangle=\left\langle X, Y \right\rangle, \forall X,Y\in\mathfrak{m}, k\in K.$$

The isotropy representation of $K$ leaves $\m$ invariant, i.e, $\Ad(K)\mathfrak{m}\subset \mathfrak{m}$ and decomposes it into irreducible submodules $$\m=\m_1\oplus\m_2\ldots\oplus\m_n,$$ and each $\m_i$ is called {\it isotropic summand}. 

If the isotropic summands of the decomposition $\m=\m_1\oplus\m_2\ldots\oplus\m_n,$ are mutually inequivalents, it follows that a $G$-invariant metric $g$ on $G/K$ is represented by a inner product $$g_{1\cdot K}=t_1Q|_{\m_1}+t_2Q|_{\m_2}+\ldots +t_nQ|_{\m_n}$$on $\m$, with $t_i$ positive constants and $Q$ is (the extension of) a inner product on $\m$, $\Ad(K)$-invariant. 

In particular, if $Q=(-B)$, with $B$ Cartan-Killing form of $G$ and $t_i=1$ for all $1\leq i\leq n$ above, the $G$-invariant metric $g$ on $G/K$ represented by the inner product $$g_{1\cdot K}=(-B)|_{\m_1}+(-B)|_{\m_2}+\ldots +(-B)|_{\m_n}$$on $\m$ is called {\it normal metric}.

\begin{remark} Let $W_{k,l}^7=G/K=SU(3)/SO(2)=SU(3)/S^1_{k,l}$ be an Aloff-Wallach space as we defined above, for some $\mathbb{Z}\ni k,l\geq 0, \,\gcd(k,l)=1$. If $\mathfrak{k}=\text{Lie}(S^1_{k,l})$, we have the decomposition $$\mathfrak{su}(3)=\mathfrak{k}\oplus\m,$$ where $\mathfrak{k}=[Z]_{\R}$ and $\m=\m_1\oplus\m_2\oplus\m_3\oplus\m_4$ is the isotropy representation of $K=S^1_{k,l}$, with isotropic summands given by $\m_1=[X_1,X_2]_{\R}$, $\m_2=[X_3,X_4]_{\R}$, $\m_3=[X_5,X_6]_{\R}$ and $\m_4=[X_0]_{\R}$. As was shown in \cite{nikoronov}, the irreducible $\Ad(K)$-modules $\m_i$ in the above decomposition are pairwise non-isomorphic for the space $W_{k,l}^7$ when $|k|\neq |l|\neq |m|\neq |k|$. The spaces where there are some pairwise isomorphic $\Ad(K)$-modules $\m_i$ are spacial cases represented by $W_{1,0}^7$ and $W_{1,1}^7$, which we will apply changes of basis in order to guarantee decompositions of their respective isotropy representation into sum of non-equivalent $\Ad(K)$-modules. In fact, $W_{1,0}$ and $W_{1,1}$ are special in the following sense: all other $W_{k,l}$ has a 4-parameter family of invariant metrics. But on $W_{1,0}$ and $W_{1,1}$ the set of invariant metrics depends on 6 and 10 parameters respectively, see details in \cite{nikoronov}.  
\end{remark} 


\section{Laplacian and Canonical Variation of Riemannian Submersions With Totally Geodesic Fibers}
\label{cvariation}


We denote by $\Delta_{g}=-\textrm{div}_g\circ\textrm{grad}_g$ the Laplacian operator of $(M,g)$ acting on $C^{\infty}(M)$. The operator $\Delta_{g}$, densely defined on $L^2(M,\textrm{vol}_g)$, is symmetric (hence closable), non-negative and has compact resolvent. Furthermore, it is well-known that $\Delta_g$ is essentially self-adjoint with this domain. We denote its unique self-adjoint extension also by $\Delta_g$. Analogously, let $\Delta_k$ be the unique self-adjoint extension of the Laplacian of the fiber $(F,k)$, where $k$ is the metric induced by $(M,g)$ on $F$.

\begin{definition} [\cite{bergery}] The \emph{vertical Laplacian} $\Delta_v$ acting on $L^2(M,\emph{vol}_g)$ is the operator defined at $p\in M$ by $$(\Delta_v\psi)(p)=(\Delta_k\psi|_{F_p})(p),$$ and the \emph{horizontal Laplacian} $\Delta_h$, acting on the same space, is defined by the difference $$\Delta_h=\Delta_g-\Delta_v.$$ 
\end{definition}

Both $\Delta_h$ and $\Delta_v$ are non-negative self-adjoint unbounded operators on $L^2(M,\textrm{vol}_g)$, but in general, are not elliptic (unless $\pi$ is a covering). We now consider the spectrum of such operators. 

As remarked above, $\Delta_g$ is non-negative and has compact resolvent, that is, its spectrum is non-negative and discrete. Since the fibers are isometrics, $\Delta_v$ also has discrete spectrum equal to the fibers. However, the spectrum of $\Delta_h$ need not be discrete.



\subsection{Canonical Variation}\label{Section4}


Let us recall the definition of {\it the canonical variation}. It will be fundamental in the determination of bifurcation and local rigid instants in our work.

\begin{definition}[\cite{bergery}]\label{variation} Let $F\cdots (M,g)\stackrel{\pi}{\rightarrow} B$ be a Riemannian submersion with totally geodesic fibers. Consider the 1-parameter family of Riemannian submersions given by \linebreak $\left\{F\cdots (M,g_t)\stackrel{\pi}{\rightarrow} B, t>0\right\}$, where $g_t\in {\mathcal{R}}^k(M)$ is defined by 
$$g_t(v,w)=\left\{
  \begin{array}{rcl} 
  t^2g(v,w),&\mbox{se} &v,w \ \ \mbox{are verticals}\\
  0,&\mbox{se} &v \ \ \mbox{is vertical and} \ \ w \ \ \mbox{is horizontal}\\
  g(v,w),&\mbox{se} &v,w \ \ \mbox{are horizontals}.\\
  \end{array}
  \right.$$

Such family of Riemannian submersions is called the \emph{canonical variation} of \linebreak$F\cdots (M,g)\stackrel{\pi}{\rightarrow} B$ or, for simplicity, we may also refer to the family of \emph{total spaces} of these submersions, i.e., the Riemannian manifolds $(M,g_t)$, as the canonical variation of $(M,g)$.
	\end{definition}
	
\begin{proposition} [\cite{bergery}] The family $\left\{F\cdots (M,g_t)\stackrel{\pi}{\rightarrow} B, t>0\right\}$ of Riemannian submersions has totally geodesic fibers, for each $t>0$. Furthermore, its fibers are isometrics to $(F,t^2k)$, where $(F,k)$ is the original fiber of $\pi:M\longrightarrow B$. 
\end{proposition}

\begin{remark}
\emph{Note that, for $a\neq b$, the metrics $g_a$ e $g_b$ are not conformal. Furthermore, for each $t>0$, $g_t$ is the unique Riemannian metric that satisfy the conditions of Definition \ref{variation}.}

\end{remark}

The following result shows how to decompose $\Delta_t$ in terms of the vertical and horizontal Laplacians.

\begin{proposition} [\cite{bettiol}] Let $\Delta_t$ the Laplacian of $(M,g_t)$. Then 

\begin{equation}
\Delta_t=\Delta_h+\frac{1}{t^2}\Delta_v=\Delta_g+(\frac{1}{t^2}-1)\Delta_v. \label{deltat}
\end{equation}
\end{proposition}

\begin{corollary} [\cite{bettiol}] \label{deltat}For each $t>0$, the following inclusion holds 
\begin{equation}
\sigma(\Delta_t)\subset\sigma(\Delta_g)+(\frac{1}{t^2}-1)\sigma(\Delta_v), \label{inclusao}
\end{equation}
where $\sigma(\Delta_t)$, $\sigma(\Delta_g)$ and $\sigma(\Delta_v)$ are the respective spectrum of $\Delta_t$, $\Delta_g$ e $\Delta_v.$ Since the above spectra are discrete, this means that every eigenvalue $\lambda(t)$ of $\Delta_t$ is of the form
\begin{equation}
\lambda^{r,j}(t)=\mu_r+(\frac{1}{t^2}-1)\phi_j, \label{kj}
\end{equation}
for some $\mu_r\in\sigma(\Delta_g)$ and $\phi_j\in\sigma(\Delta_v).$ 
\end{corollary}

\begin{corollary} \label{primeiroaut1}If $\lambda_1(t)$ is the first positive eigenvalue of $\Delta_t$, then $$\lambda_1(t)\geq \mu_1, \forall \ \ 0<t\leq 1,$$ where $\mu_1$ is the first positive eigenvalue of $\Delta_g$.
\end{corollary}
{\bf Proof:} Since $\lambda_1(t)=\mu_r+(\frac{1}{t^2}-1)\phi_j$, for some $\mu_r\in \sigma(\Delta_g)$ and $\phi_j\in \sigma(\Delta_v)$, if $0<t\leq 1$, $(\frac{1}{t^2}-1)\phi_j\geq 0$ and $\lambda_1(t)=\mu_r+(\frac{1}{t^2}-1)\phi_j\geq\mu_r\geq\mu_1,$ where $\mu_1\in\sigma(\Delta_g)$ is the first positive eigenvalue of the operator $\Delta_g$. \cqd

We also remark that not all possible combinations of $\mu_r$ and $\phi_j$ in the expression \eqref{kj} give rise to an eigenvalue of $\Delta_t$. In fact, this only happens when the total space $(M,g)$ of the submersion is a Riemannian product. To determine which combinations are allowed in general is difficulty and depends on the global geometry of the submersion. Moreover,the ordering of the eigenvalues of $\Delta_t$ may change with  $t$. 

We have the following useful properties of the spectrum of the Laplacian operator $\Delta_t$ of $g_t$.

\begin{proposition} [\cite{pacificjournal}] Using the same notations above, one has that $$\sigma(\Delta_B)\subset\sigma(\Delta_t),$$ for all $t>0.$
\end{proposition}


We remark that the spectrum of $\Delta_h$, $\sigma(\Delta_h)$, contains but not coincides with the spectrum of the basis $B$. In fact, if $\overline{f}$ is a $C^{\infty}$ function on the basis $B$, then $$(\Delta_{g'}\overline{f})\circ\pi=\Delta_g(\overline{f}\circ\pi)=\Delta_h(\overline{f}\circ\pi),$$ where $\Delta_{g'}$ is the Laplacian operator acting on functions in $C^{\infty}(B,g')$.

\begin{corollary} \label{primeiroaut2}Denoting by $\beta_1$ the first positive eigenvalue of $\Delta_B$ and by $\lambda_1(t)$ the first positive eigenvalue of $\Delta_t$, the following inequality holds $$\lambda_1(t)\leq \beta_1, \ \ \forall \, t>0.$$
\end{corollary}
  
As consequence of the Corollaries \ref{primeiroaut1} and \ref{primeiroaut2}, we have that $\mu_1\leq\lambda_1(t)\leq \beta_1$, where $\mu_1$ the first positive eigenvalue of the Laplacian $\Delta_g$ on the total space $M$. 

When $j=0$ in the expression \eqref{kj}, if $\lambda^{r,0}(t)=\mu_r\in\sigma(\Delta_g)$ remains an eigenvalue of $\Delta_t$ for $t\neq 1$, such eigenvalues will be called {\it constant eigenvalue of $\Delta_t$}, since they are independent of $t$. We stress that $\lambda^{r,0}(t)$ is not necessarily a constant eigenvalue of $\Delta_t$ for all $k$. A simple criterion to determinate when $\lambda^{r,0}(t)\in\sigma(\Delta_t)$ is given in the following:

\begin{proposition} [\cite{pacificjournal}] \label{constanteigen}$\mu_r=\lambda^{r,0}(t)\in\sigma(\Delta_t)$ for $t\neq 1$ $\Leftrightarrow$ $\mu_r\in\sigma(\Delta_B)$.
\end{proposition}

\subsubsection{Homogeneous Fibration}\label{Section2.2}
In our main result, the canonical variations $g_t$, obtained from a Riemannian submersion with totally geodesic fibers, are {\it homogeneous metrics}, which are trivial solutions to the Yamabe problem, since every homogeneous metric has constant scalar curvature. This makes these metrics good candidates for admitting other solutions in their conformal class. 

The homogeneous fibrations are obtained from the following construction. Let $K\subsetneq H\subsetneq G$ be compact connected Lie groups, such that $\dim K/H\geq 2$. Consider the natural fibration
$$
\begin{array}{cccc} 
\pi :&\! G/K &\!\longrightarrow &\!G/H\\
     &\! \alpha K &\!\mapsto &\! \alpha H,
		\end{array}
		$$
with fibers $H/K$ $(\textrm{em} \ \ eK)$ and structural group $H$. More precisely, $\pi$ is the associated bundle with fiber $H/K$ to the $H$-principal bundle $p:G\longrightarrow G/H$. 	

Let $\lie$ be the Lie algebra of $G$ and $\mathfrak{h}\supset\mathfrak{k}$ the Lie algebras of $H$ and $K$, respectively. Given a inner product on $\lie$, determined by the Cartan-Killing form $B$ of $\lie$, since $K,H$ and $G$ are compacts, we can consider a $\Ad_G(H)$-invariant orthogonal complement $\mathfrak{q}$ to $\mathfrak{h}$ in $\lie$, i.e., $[\h,\mathfrak{q}]\subset\mathfrak{q}$, and a $\Ad_G(K)$-invariant orthogonal complement $\mathfrak{p}$ to $\sub$ in $\h$, i.e., $[\sub,\mathfrak{p}]\subset\mathfrak{p}$. It follows that $\mathfrak{p}\oplus\mathfrak{q}$ is a $\Ad_G(K)$-invariant orthogonal complement to $\sub$ in $\lie$. 

The $\Ad_G(H)$-invariant inner product $(-B)|_{\mathfrak{q}}$ on $\mathfrak{q}$ define a $G$-invariant Riemannian metric $\breve{g}$ on $G/H$, and the inner product $(-B)|_{\mathfrak{p}}$, $\Ad_G(K)$-invariant on $\mathfrak{p}$, define on $H/K$ a $H$-invariant Riemannian metric $\hat{g}$ on $H/K$. Finally, the orthogonal direct sum of these inner products on $\mathfrak{p}\oplus\mathfrak{q}$ define a $G$-invariant metric $g$ on $G/K$, determined by 
\begin{equation}
g(X+V,Y+W)_{eK}=(-B)|_{\mathfrak{q}}(X,Y) + (-B)|_{\mathfrak{p}}(V,W), \label{metric} 
\end{equation}
for all $X,Y\in\mathfrak{q}$ and $V,W\in\mathfrak{p}$; $g$ is a normal homogeneous metric on $G/K$ and the $(-B)$-orthogonal direct sum $\m=\p\oplus\q$ is the isotropy representation of $K$.

\begin{theorem}[(\cite{besse}, p. 257)] The map $\pi:(G/K,g)\longrightarrow (G/H,\breve{g})$ is a Riemannian submersion with totally geodesic fibers and isometric to $(H/K,\hat{g})$.
\end{theorem} 

If we take for each $t>0$ the metric $g_t$ that corresponds to the inner product on $\mathfrak{p}\oplus\mathfrak{q}$
\begin{equation}
 \left\langle \cdot,\cdot\right\rangle_t=(-B)|_{\mathfrak{q}}+(-t^2B)|_{\mathfrak{p}}, \label{gt} 
\end{equation}
where $B$ is the Cartan-Killing form of $\lie$, we have that the map $\pi_t:(G/K,g_t)\longrightarrow (G/H,\breve{g})$ is a Riemannian submersion with totally geodesic fibers and isometric to $H/K$ provided with the induced metric $-t^2B|_{\mathfrak{p}}$. Hence, we obtain the canonical variation of the original homogeneous fibration $\pi:(G/K,g)\longrightarrow (G/H,\breve{g})$.

\subsection{{Canonical Variations of Normal Metrics on Aloff-Wallach Spaces}}\label{Section2.3}

Remembering that the Aloff-Wallach spaces $W_{k,l}^7=SU(3)/SO(2)$ are defined by an embedding of the circle $SO(2)=S^1$ into $SU(3)$ of the type $$i_{l,k}:e^{2\pi il\theta}\mapsto \textrm{diag}(e^{2\pi ik\theta},e^{2\pi il\theta},e^{2\pi im\theta}),$$ where $k,l$ are integers with greatest common divisor 1 and $k+l+m=0$. By using Weil group of $SU(3)$ one can assume that $k\geq l\geq 0$ (see \cite{nikoronov}).


\subsubsection{$(W_{k,l}^7,\cvg), \gcd(k,l)=1, k>l>0$:}

We observe that for $\gcd(k,l)=1, k>l>0$ and considering the Aloff-Wallach spaces $W_{k,l}^7=G/K=SU(3)/S_{k,l}^1$, if $\mathfrak{k}=\text{Lie}(S^1_{k,l})$, we have the $(-B)$-orthonormal decomposition, $B$ Cartan-Killing form of $\mathfrak{su}(3)$ given by $B(X,Y)=6\text{ReTr(XY)}$, $$\mathfrak{su}(3)=\mathfrak{k}\oplus\m,$$ where $\mathfrak{k}=[Z]_{\R}$ and $\m=\m_1\oplus\m_2\oplus\m_3\oplus\m_4$ is the isotropy representation of $K=S^1_{k,l}$, with isotropic summands given by $\m_1=[X_1,X_2]_{\R}$, $\m_2=[X_3,X_4]_{\R}$, $\m_3=[X_5,X_6]_{\R}$ and $\m_4=[X_0]_{\R}$, so that theses $\Ad(K)$-modules are pairwise non-equivalent. 
 
Setting $H=U(2)$ and $K=S^1_{k,l}$, we have $K\subsetneq H\subsetneq G=SU(3)$, with $G,H$ and $K$ compact connected Lie groups. Consider the canonical map 
\begin{center}
$U(2)/S^1_{k,l}\cdots W^7_{k,l}=SU(3)/S^1_{k,l} \rightarrow SU(3)/U(2)=\mathbb{CP}^{2}.$
\end{center}  

Let $\lie =\mathfrak{su}(3)$, $\h=\mathfrak{su}(2)$ and $\sub$ be the Lie algebras of $G$, $H$ and $K$, respectively. As we saw previously in Section \ref{Section2.2}, since $K,H$ and $G$ are compacts, we can  consider a $\Ad(H)$-invariant $(-B)$-orthogonal complement $\mathfrak{q}$ to $\h$ in $\lie$, and a $\Ad(K)$-invariant $(-B)$-orthogonal complement $\mathfrak{p}$ to $\sub$ in $\h$.

The Lie algebra $\lie$ decomposes into the sum $$\lie=\mathfrak{su}(3)=\sub\oplus \mathfrak{p}\oplus\mathfrak{q}=\sub\oplus\m,$$ hence
\begin{equation}
\pi:(W_{k,l}^7,g)\longrightarrow (SU(3)/U(2)=\mathbb{CP}^{2},\breve{g}) \label{homfibw}
\end{equation}
is a Riemannian submersion with totally geodesic fibers isometric to $(U(2)/S^1_{k,l},\hat{g})$, with $g$ the normal metric determined by the inner product $(-B)|_{\m}$, $\hat{g}$ the metric given by $(-B)|_{\mathfrak{p}}$ and $\breve{g}$ defined by the inner product $(-B)|_{\mathfrak{q}}$.

The fiber $U(2)/S_{k,l}^1=S^3/\text{diag}(z^k,z^l)=S^3/\mathbb{Z}_{k+l}, \; z^{k+l}=1$, of the original fibration is a 3-dimensional lens space and from the brackets $[X_i,X_j]$ above, we can conclude that the vertical component, tangente to the fiber, is given by $\mathfrak{p}=\m_3\oplus\m_4.$ 
The Riemannian metric $\hat{g}$ on the fiber $U(2)/S_{k,l}^1$, represented by $(-B)|_{\mathfrak{p}}$, is the normal metric represented by the inner product $$g_{eK}=(-B)|_{\m_3}+(-B)|_{\m_4}$$and the homogeneous metric $\breve{g}$ defined by the inner product $(-B)|_{\mathfrak{q}}$, $\mathfrak{q}=\m_1\oplus\m_2$, makes the basis $(SU(3)/U(2)=\mathbb{CP}^{2},\breve{g})$ an isotropy irreducible compact symmetric space.

Hence, by scaling the normal metric $$g_{eK}=(-B)|_{\m_1}+(-B)|_{\m_2}+(-B)|_{\m_3}+(-B)|_{\m_4}$$ on the total space $W_{k,l}^7$ in the direction of the fibers by $t^2$, i.e., multiplying by $t^2, t>0$, the parcels $(-B)|_{\m_3}+(-B)|_{\m_4}$ in the expression of $g$ we obtain the canonical variation $\cvg$ of the metric $g$, according \eqref{gt} 

$${(\cvg)}_{eK}=(-B)|_{\m_1}+(-B)|_{\m_2}+t^2(-B)|_{\m_3}+t^2(-B)|_{\m_4}$$

\subsubsection{The special case $(W_{1,0}^7,\cvh)$: \label{10}}

If $(k,l,m)=(1,0,-1)$, we have the decomposition $$\mathfrak{su}(3)=\mathfrak{k}\oplus\m,$$ where $\mathfrak{k}=[Z]_{\R}$ and $\m=\m_1\oplus\m_2\oplus\m_3\oplus\m_4$ is the isotropy representation of $K=S^1_{k,l}$, with isotropic summands given by $\m_1=[X_1,X_2]_{\R}$, $\m_2=[X_3,X_4]_{\R}$, $\m_3=[X_5,X_6]_{\R}$ and $\m_4=[X_0]_{\R}$, so that $\m_1$ and $\m_3$ are isomorphic to each other. In fact, every isomorphism $\varphi:\m_1\longrightarrow\m_3$ has the form (see \cite{nikoronov}) $$\varphi(aX_1+bX_2)=a(\alpha X_5+\beta X_6)+b(-\beta X_5+\alpha X_6).$$
 
Setting $H=U(2)$ and $K=S^1_{1,0}$, we have $K\subsetneq H\subsetneq G=SU(3)$, with $G,H$ and $K$ compact connected Lie groups. Consider the canonical map 
\begin{center}
$U(2)/S^1_{1,0}=S^3\cdots W^7_{1,0}=SU(3)/S^1_{1,0} \rightarrow SU(3)/U(2)=\mathbb{CP}^{2}.$
\end{center}  


Let $\lie = \mathfrak{su}(3)$, $\h = \mathfrak{su}(2)$, and $\sub$ be the Lie algebras of $G$, $H$, and $K$, respectively. As we saw previously in Section \ref{Section2.2}, since $K$, $H$, and $G$ are compact, we can consider an $\Ad(H)$-invariant $(-B)$-orthogonal complement $\mathfrak{q}$ to $\h$ in $\lie$, and an $\Ad(K)$-invariant $(-B)$-orthogonal complement $\mathfrak{p}$ to $\sub$ in $\h$.

The Lie algebra $\lie$ decomposes into the sum $$\lie=\mathfrak{su}(3)=\sub\oplus \mathfrak{p}\oplus\mathfrak{q}=\sub\oplus\m,$$ hence
\begin{equation}
\pi:(W_{1,0}^7,g)\longrightarrow (SU(3)/U(2)=\mathbb{CP}^{2},\breve{g}) \label{homfibw}
\end{equation}
is a Riemannian submersion with totally geodesic fibers isometric to $(U(2)/S^1_{1,0},\hat{g})$, with $g$ the normal metric determined by the inner product $(-B)|_{\m}$, $\hat{g}$ the metric given by $(-B)|_{\mathfrak{p}}$ and $\breve{g}$ defined by the inner product $(-B)|_{\mathfrak{q}}$. The fiber $U(2)/S_{1,0}^1$ of the original fibration is the 3-dimensional round sphere $S^3$.

In order to determine the vertical component $\mathfrak{p}$, tangent to the fiber $U(2)/S_{1,0}^1=S^3$, as well as the horizontal one $\mathfrak{q}$ tangent to the basis $SU(3)/U(2)=\mathbb{CP}^{2}$, it is more appropriate the following description of the invariant metrics on the total space $W_{1,0}^7$. Namely, for $X,Y\in\m$, an inner product on $\m$ is given by $\langle X,Y\rangle=(A(X),Y)$, where $A$ is a positive-definite operator in $S^2(\m)^{\Ad (H)}$, the space of $\Ad (H)$-equivariant, $(\cdot, \cdot)$-self-adjoint endomorphisms of $\m$, with $\det A=1$. Therefore, applying the action of the gauge group of the set of $\Ad (H)$-invariant inner products on $\m$, it is possible to find, for a given invariant metric, another one which is more convenient, such as when $(\cdot,\cdot)$ and $\langle \cdot,\cdot\rangle$ have diagonal forms, enabling us to find irreducible modules $\m_i$ non-equivalent each other in the new decomposition of $\m$. For a more detailed description of this process, see \cite{nikoronov}. 

Actually, we can obtain the matrix $A$ given by 

$$A=\left( 
\begin{array}{cccc}
D_1 & 0 & 0 & 0 \\
0 & D_2 & 0 & C \\
0 & 0 & D_3 & 0 \\
0 & C^T & 0 & D_4 \\
\end{array}
\right), 
$$ 
with respect to  the basis $\{X_0, X_1, X_2, X_3, X_4, X_5, X_6\}$ such that $\langle X,Y\rangle=(A(X),Y)$, where $(X,Y)=-\frac{1}{2}\text{ReTr}(XY),$ $$D_1=\text{diag}(z_1), \quad D_2=\text{diag}(z_2,z_2), \quad D_3=\text{diag}(z_3,z_3), $$ $$ D_4=\text{diag}(z_4,z_4), \quad 
C= \left( \begin{array}{cc}
\beta & \gamma \\
-\gamma & \beta \\
\end{array}
\right).
$$

If we denote by $N_G(H)$ the normalizer of $H=SO(2)$ in $G=SU(3)$, we have that the gauge group of the $\Ad (H)$-invariant inner products on $\m$ is isomorphic to $SO(2)$. By considering the natural action of this group on $W_{1,0}^7=SU(3)/SO(2)$, given by $(gH,nH)\rightarrow gnH$, it can be determined the basis $\{X_0,Y_1,Y_2,X_3,X_4,Y_3,Y_4\}$, where $$Y_1=\cos(\alpha)X_1 + \sin(\alpha)X_5, \quad Y_2=\cos(\alpha)X_2 + \sin(\alpha)X_6,$$ $$Y_3=-\sin(\alpha)X_1 + \cos(\alpha)X_5, \quad Y_4=-\sin(\alpha)X_2 + \cos(\alpha)X_6,$$ in which the matrix $A$ of the metric $\langle \cdot,\cdot\rangle=(A \cdot,\cdot)$ makes both the forms $\langle \cdot,\cdot\rangle$ and $(\cdot,\cdot)$ diagonal, for some $\alpha\in \R$. 

If we take $\mathfrak{p}_1=[Y_1,Y_2]_{\R}$, $\mathfrak{p}_2=[Y_3]_{\R}$, $\mathfrak{p}_3=\m_2$, $\mathfrak{p}_4=\m_4$, one has that each $\mathfrak{p}_i$ is a $\Ad(H)$-invariant irreducible module and theses submodules are inequivalent each other. It follows that the set of invariant metrics on $W_{1,0}^7$ has the form $$g_{1\cdot K}=t_1Q|_{\mathfrak{p}_1}+t_2Q|_{\mathfrak{p}_2}+t_3Q|_{\mathfrak{p}_3}+t_4Q|_{\mathfrak{p}_4},$$ where $t_i$, $1\leq i\leq 4$, are positive numbers and $Q(X,Y)=(X,Y)=-\frac{1}{2}\text{ReTr}(XY)$. In particular, the normal metric on $W_{k,l}^7$ is given by
$$g_{1\cdot K}=Q|_{\mathfrak{p}_1}+Q|_{\mathfrak{p}_2}+Q|_{\mathfrak{p}_3}+Q|_{\mathfrak{p}_4}.$$

The fiber $U(2)/S_{1,0}^1=S^3$ of the original fibration is a 3-dimensional round sphere and from the brackets $[X_i,X_j]$ above, we can conclude that the vertical component, tangente to the fiber, is given by $\mathfrak{p}=\mathfrak{p}_3\oplus\mathfrak{p}_4.$ 
The Riemannian metric $\hat{g}$ on the fiber $S^3$, represented by $(-B)|_{\mathfrak{p}}$, is the normal metric represented by the inner product $$g_{eK}=(-B)|_{\mathfrak{p}_3}+(-B)|_{\mathfrak{p}_3}$$and the homogeneous metric $\breve{g}$ defined by the inner product $(-B)|_{\mathfrak{q}}$, $\mathfrak{q}=\mathfrak{p}_1\oplus\mathfrak{p}_2$, makes the basis $(SU(3)/U(2)=\mathbb{CP}^{2},\breve{g})$ an isotropy irreducible compact symmetric space.

Hence, by scaling the normal metric $$g_{eK}=(-B)|_{\mathfrak{p}_1}+(-B)|_{\mathfrak{p}_2}+(-B)|_{\mathfrak{p}_3}+(-B)|_{\mathfrak{p}_4}$$ on the total space $W_{k,l}^7$ in the direction of the fibers by $t^2$, i.e., multiplying by $t^2, t>0$, the parcels $(-B)|_{\mathfrak{p}_3}+(-B)|_{\mathfrak{p}_4}$ in the expression of $g$ we obtain the canonical variation $\cvh$ of the metric $g$, according \eqref{gt}.

$${(\cvh)}_{eK}=(-B)|_{\mathfrak{p}_1}+(-B)|_{\mathfrak{p}_2}+t^2(-B)|_{\mathfrak{p}_3}+t^2(-B)|_{\mathfrak{p}_4}.$$

\subsubsection{The special case $(W_{1,1}^7,\cvk)$:\label{11}} 

If $(k,l,m)=(1,1,-2)$, we have the decomposition $$\mathfrak{su}(3)=\mathfrak{k}\oplus\m,$$ where $\mathfrak{k}=[Z]_{\R}$ and $\m=\m_1\oplus\m_2\oplus\m_3\oplus\m_4$ is the isotropy representation of $K=S^1_{1,1}$, with isotropic summands given by $\m_1=[X_1,X_2]_{\R}$, $\m_2=[X_3,X_4]_{\R}$, $\m_3=[X_5,X_6]_{\R}$ and $\m_4=[X_0]_{\R}$, so that $\m_2$ and $\m_3$ are isomorphic to each other. In fact, every isomorphism $\varphi:\m_2\longrightarrow\m_3$ has the form (see \cite{nikoronov}) $$\varphi(aX_3+bX_4)=a(\alpha X_5+\beta X_6)+b(\beta X_5-\alpha X_6).$$
 
Setting $H=U(2)$ and $K=S^1_{1,1}$, we have $K\subsetneq H\subsetneq G=SU(3)$, with $G,H$ and $K$ compact connected Lie groups. Consider the canonical map 
\begin{center}
$U(2)/S^1_{1,1}=U(2)/Z(U(2))=SO(3)\cdots W^7_{1,1}=SU(3)/S^1_{1,1} \rightarrow SU(3)/U(2)=\mathbb{CP}^{2}.$
\end{center}  

Let $\lie =\mathfrak{su}(3)$, $\h=\mathfrak{su}(2)$ and $\sub$ be the Lie algebras of $G$, $H$ and $K$, respectively. As we saw previously in \ref{Section2.2}, since $K,H$ and $G$ are compacts, we can consider a $\Ad(H)$-invariant $(-B)$-orthogonal complement $\mathfrak{q}$ to $\h$ in $\lie$, and a $\Ad(K)$-invariant $(-B)$-orthogonal complement $\mathfrak{p}$ to $\sub$ in $\h$.

The Lie algebra $\lie$ decomposes into the sum $$\lie=\mathfrak{su}(3)=\sub\oplus \mathfrak{p}\oplus\mathfrak{q}=\sub\oplus\m,$$ hence
\begin{equation}
\pi:(W_{1,1}^7,g)\longrightarrow (SU(3)/U(2)=\mathbb{CP}^{2},\breve{g}) \label{homfibw}
\end{equation}
is a Riemannian submersion with totally geodesic fibers isometric to $(U(2)/S^1_{1,1},\hat{g})=(SO(3),\hat{g})$, with $g$ the normal metric determined by the inner product $(-B)|_{\m}$, $\hat{g}$ the metric given by $(-B)|_{\mathfrak{p}}$ and $\breve{g}$ defined by the inner product $(-B)|_{\mathfrak{q}}$.

As in the previous case, we must determine another decomposition for $\m$, in such a way that their isotropic summands are not pairwise equivalent to each other. For $X,Y\in\m$, remember that an inner product on $\m$ is given by $\langle X,Y\rangle=(A(X),Y)$, where $A$ is a positive-definite operator in $S^2(\m)^{\Ad (H)}$, the space of $\Ad (H)$-equivariant, $(\cdot, \cdot)$-self-adjoint endomorphisms of $\m$, with $\det A=1$. If we denote by $N_G(H)$ the normalizer of $H=U(2)$ in $G=SU(3)$, we have that in this case the gauge group $N_G(H)/H$ of the $\Ad (H)$-invariant inner products on $\m$ is isomorphic to $U(2)$. Applying the action of the gauge group of the set of $\Ad (H)$-invariant inner products on $\m$, in \cite{nikoronov} was determined, the following decomposition of $\m$ with five isotropic summands, pairwise inequivalent to each other. Actually, in \cite{nikoronov} one has the basis $\{X_3,X_4,X_5,X_6,Y_1,Y_2,Y_3\}$, with $Y_1,Y_2,Y_3$ obtained from the matrix transformation from the basis $\{X_1,X_2,X_3\}$ to the basis $\{Y_1,Y_2,Y_3\}$, given by      

$$T=\left( 
\begin{array}{ccc}
-\dfrac{\sin\varphi\sin\psi}{\rho} & \dfrac{\cos\varphi\cos\psi}{\rho} & \dfrac{\cos\varphi\sin\psi}{\rho} \\
-\dfrac{\sin\varphi\cos\varphi\cos\psi}{\rho} & -\dfrac{\sin\psi}{\rho} & \dfrac{\cos^2\varphi\cos\psi}{\rho} \\
\cos\varphi & 0 & \sin\varphi \\
\end{array}
\right), 
$$ 
where $\rho=\sqrt{\cos^2\varphi+\sin^2\varphi\sin^2\psi}$, $0\leq \varphi\leq\arccos (\frac{\sqrt{3}}{3}),$ the angle between $Y_3$ and $X_0$, $0\leq\psi\leq \frac{\pi}{4}$ the angle between $X_1$ and the projection of $Y_1$ on the plane determined by $X_1$ and $X_2$.

With respect to the basis $\{X_3,X_4,X_5,X_6,Y_1,Y_2,Y_3\}$, the matrix $A$ of the metric $\langle \cdot,\cdot\rangle=(A \cdot,\cdot)$ makes both the forms $\langle \cdot,\cdot\rangle$ and $(\cdot,\cdot)$ diagonal.

If we take $\mathfrak{p}_1=[X_3,X_4]_{\R}$, $\mathfrak{p}_2=[X_5,X_6]_{\R}$, $\mathfrak{p}_3=[Y_1]_{\R}$, $\mathfrak{p}_4=[Y_2]_{\R}$,  $\mathfrak{p}_5=[Y_3]_{\R}$ one has that each $\mathfrak{p}_i$ is a $\Ad(H)$-invariant irreducible module and theses submodules are inequivalent each other. It follows that the set of invariant metrics on $W_{1,1}^7$ has the form $$g_{1\cdot K}=t_1Q|_{\mathfrak{p}_1}+t_2Q|_{\mathfrak{p}_2}+t_3Q|_{\mathfrak{p}_3}+t_4Q|_{\mathfrak{p}_4}+t_5Q|_{\mathfrak{p}_5},$$ where $t_i$, $1\leq i\leq 5$, are positive numbers and $Q(X,Y)=(X,Y)=-\frac{1}{2}\text{ReTr}(XY)$. In particular, the normal metric on $W_{1,1}^7$ is given by
$$g_{1\cdot K}=Q|_{\mathfrak{p}_1}+Q|_{\mathfrak{p}_2}+Q|_{\mathfrak{p}_3}+Q|_{\mathfrak{p}_4}+Q|_{\mathfrak{p}_5}.$$

On the fiber $U(2)/S_{1,1}^1=SO(3)$ of the original fibration, from the brackets $[X_i,X_j]$ above, we can conclude that the vertical component, tangente to the fiber, is given by $\mathfrak{p}=\mathfrak{p}_3\oplus\mathfrak{p}_4\oplus\mathfrak{p}_5.$ 
The Riemannian metric $\hat{g}$ on the fiber $U(2)/S_{1,1}^1=SO(3)$, represented by $(-B)|_{\mathfrak{p}}$, is the normal metric represented by the inner product $$g_{eK}=(-B)|_{\mathfrak{p}_3}+(-B)|_{\mathfrak{p}_4}+(-B)|_{\mathfrak{p}_5}$$and the homogeneous metric $\breve{g}$ defined by the inner product $(-B)|_{\mathfrak{q}}$, $\mathfrak{q}=\mathfrak{p}_1\oplus\mathfrak{p}_2$, makes the basis $(SU(3)/U(2)=\mathbb{CP}^{2},\breve{g})$ an isotropy irreducible compact symmetric space.

Hence, by scaling the normal metric $$g_{eK}=(-B)|_{\mathfrak{p}_1}+(-B)|_{\mathfrak{p}_2}+(-B)|_{\mathfrak{p}_3}+(-B)|_{\mathfrak{p}_4}+(-B)|_{\mathfrak{p}_5}$$ on the total space $W_{1,1^7}$ in the direction of the fibers by $t^2$, i.e., multiplying by $t^2, t>0$, the parcels $(-B)|_{\mathfrak{p}_3}+(-B)|_{\mathfrak{p}_4}+(-B)|_{\mathfrak{p}_5}$ in the expression of $g$ we obtain the canonical variation $\cvh$ of the metric $g$, according \eqref{gt}.

$${(\cvk)}_{eK}=(-B)|_{\mathfrak{p}_1}+(-B)|_{\mathfrak{p}_2}+t^2(-B)|_{\mathfrak{p}_3}+t^2(-B)|_{\mathfrak{p}_4}+t^2(-B)|_{\mathfrak{p}_5}.$$

\subsection{Spectra of Aloff-Wallach spaces} \label{spectraspaces}

In \cite{Urakawa84}, H. Urakawa describes the spectrum of the Laplacian defined on $C^{\infty}(W_{k,l}^7,g)$, with $W_{k,l}^7=SU(3)/S^1_{k,l}$ equipped with a normal homogeneous metric $g$ determined by the inner product $(X,Y)=-\frac{1}{2}\text{ReTr}(XY)$. In this case, we have the following description of $\sigma(\Delta_g)$.

\begin{theorem}[H. Urakawa \cite{Urakawa84}] The spectrum of Laplacian of the Riemannian manifold $(W_{k,l}^7=SU(3)/S^1_{k,l},g)$ is given as follows: 

Eigenvalues: $f(n_1,n_2)=\frac{4}{3}(m_1^2+m_2^2-m_1m_2+3m_1)$ 
where $m_1:=n_1+n_2$, $m_2=n_2$ and $n_1$ and $n_2$ run over the set of nonnegative integers satisfying $S_{n_1,n_2}^{k,l}\neq 0$, with $S_{n_1,n_2}^{k,l}$ being the number of all the integer solutions $(p',q,r)$ of the equations 

$$
\left\{
\begin{array}{cc}
 kn_1-ln_2-(2k+l)p'+(-k+l)q+(k+2l)r = 0 \\
 0\leq p'\leq n_1, \, 0\leq q\leq n_2, \, \emph{and} \; 0\leq r\leq p'+(n_2-q).  \\
 \end{array}
\right.
$$

\end{theorem}

Note that our metric $g$ in this paper is $\frac{1}{2}$ times the Riemannian metric in \cite{Urakawa84}, so that the eigenvalues of $W_{k,l}^7$ here are 2 times the ones of \cite{Urakawa84}. On the other hand, from \cite{Urakawa}, we have that the first positive eigenvalue in this case is equal to $\lambda_1=12$. Remembering that $\text{Spec}(\Delta_{\alpha g})=\frac{1}{\alpha}\text{Spec}(\Delta_{g})$, it follows that the first positive eigenvalue $\lambda_1$ of the normal metric on $W^7_{k,l}$ induced by the Cartan-Killing form $B$ of $\mathfrak{su}(3)$, given by $B(X,Y)=-6\text{ReTr(XY)}$, is $\lambda_1=1$

Let $G$ be a compact simply connected simple Lie group, $H$ closed subgroup of $G$. Let $\lie,\h$ the Lie algebras of $G$ and $H$, respectively, and $\lie=\h\oplus\q$, the Cartan decomposition. The inner product on $\q$ is $(-B)|_{\q}$, where $B$ is the Cartan-Killing form of $\lie$. Let $\breve{g}$ the $G$-invariant Riemannian metric on $G/H$ induced by $B$. Then, it is known that the spectrum of the Laplacian of $(G/H,\breve{g})$ is given by 
\begin{equation}
\mu(\Lambda)=-B(\Lambda+2\delta,\Lambda),\label{specsym}
\end{equation} with multiplicities 
\begin{equation}
d_{\Lambda}=\prod_{\alpha\in R^{+}}\dfrac{-B(\Lambda+\delta,\alpha)}{-B(\delta,\alpha)}.\label{multsym}
\end{equation}
Here, $\Lambda$ varies over the set $D(G,H)$ of the highest weights of all spherical representations of $(G,H)$, $\delta$ is equal to the sum of the positive roots $\alpha\in R^{+}$ of the complexification $\lie^{\mathbb{C}}$ of $\lie$ relative to the maximal abelian subalgebra $\sub^{\mathbb{C}}$ of $\lie^{\mathbb{C}}$ and $d_{\Lambda}$ is the dimension of the irreducible spherical representation of $(G,H)$ with highest weight $\Lambda.$
 
\subsubsection{Spectrum of the Complex Projective Space, $\sigma(\Delta_{\mathbb{CP}^2})$}
The basis space $(G/H,\breve{g})=(SU(3)/ U(2))=\mathbb{CP}^2,\breve{g})$ of the homogeneous fibration $$\pi:(W_{k,l}^7,g)\longrightarrow \mathbb{CP}^2,$$ is the homogeneous realization of the complex projective space $\mathbb{CP}^2$ equipped with a metric homothetical to the Fubini-Study one. In fact, the induced metric $(-B)|_{\q}$ on $\mathbb{CP}^2$, $B$ being the Cartan-Killing form given by $B(X,Y)=-6\text{ReTr}(XY),X,Y\in\mathfrak{su}(3)$, is simply $3$ times the Fubini-Study metric. According to \cite{taniguchi}, the spectra of the Laplacian acting on functions on the complex projective space $\mathbb{CP}^2$ equipped with the Fubini-Study metric is $$\sigma(\Delta_{FS})=\{\xi_k=k(k+2);k\in\mathbb{N}\}.$$ Hence, the spectrum of the Laplacian on $\mathbb{CP}^2=SU(3)/U(2))$ with the metric $\breve{g}$ represented by the inner product $-B|_{\q}$ is
\begin{equation}
 \sigma(\Delta_{\breve{g}})=\left\{\beta_k=\dfrac{\xi_k}{3}=\dfrac{k(k+2)}{3};k\in\mathbb{N}\right\}.\label{specAn}
\end{equation} Note that the first positive eigenvalue in this case is $\beta_1=1.$

\subsubsection{Spectrum of homogeneous Lens Spaces, $\sigma(\Delta_{S^3/\mathbb{Z}_n})$}

In \cite{Sakai}, Sakai determined the spectrum of the homogeneous lens spaces $(M=S^{2n-1}/G,g)$, with $(S^{2n-1},g_0)$ imbedded in $\mathbb{C}^n=\R^{2n}$, $g_0$ the round metric and $G$ cyclic group of order $p$, generated by some $T\in SO(2n)$ given by $$T:(z_1,z_2,\ldots,z_n)\longrightarrow (e^{2\pi/p}z_1,e^{2\pi/p}z_2, \ldots, e^{2\pi/p}z_n).$$

\begin{theorem}[\cite{Sakai}]
The eigenvalues of $\Delta_g$ are given by $\lambda_{m,r}=(2m+rp)(2n-2+2m+rp)$, $m=0,1,2,\ldots;\, r=0,1,2,\ldots$    
\end{theorem}

For $U(2)/S_{k,l}^1=S^3/\text{diag}(z^k,z^l)=S^3/\mathbb{Z}_{k+l}, \; z^{k+l}=1$, we identify $G$ with the cyclic group $\mathbb{Z}_{k+l}$ of order $k+l$. 
According Table 1 of \cite{Sakai}, the first positive eigenvalue of $\Delta_g$ must be $\lambda_1=8$.


\section{Bifurcation and Local Rigidity Instants for
Canonical Variations on Maximal Flag
Manifolds}
\label{section5}

We will determine the bifurcation and local rigidity instants for the canonical variations $\cvg,\cvh,\cvk$ defined in the Section \ref{Section2.3}. The criterion used to find such instants is based on comparison between an expression of the eigenvalues of the Laplacian relative to the respective canonical variation and a multiple of its scalar curvature. This method is related to the notion of Morse index.   

We will use an expression for the scalar curvature of each canonical variation described above by mean a general formula due M. Wang and W. Ziller in \cite{WZ} that can be applied to calculate the scalar curvature of homogeneous metrics on reductive homogeneous spaces. 

Thereafter, necessary conditions for classifications of bifurcation and local rigidity instants in the interval $]0,\infty[$ can be deduced by using the expressions of the scalar curvature in all cases.

\subsection{Scalar Curvature}
The general formula of the scalar curvature for reductive  homogeneous spaces which we use here is obtained in \cite{WZ}.

Let $G$ be a compact connected Lie group and $K\subset G$ a closed subgroup of $G$. Assume that $K$ is connected, which corresponds to the case that $G/K$ is simply connected. Let $\m$ be the $(-B)$-orthogonal complement to $\sub$ in $\lie$, where $\lie$ and $\sub$ are the Lie algebras of $G$ and $K$, respectively, and $B$ is the Cartan-Killing form of $\lie$. It is known that the isotropy representation $\m$ of $K$ decomposes into a direct sum of inequivalent irreducible submodules, $$\m=\m_1\oplus\m_2\oplus\ldots\oplus\m_r,$$ with $\m_1,\m_2,\ldots,\m_r$ such that $\Ad(K)\m_i\subset\m_i$, for all $i=1,\ldots,r$. Thus, each $G$-invariant metric on $G/K$ can be represented by a inner product on $\m$ given by  $$t_1(-B)|_{\m_1}+t_2(-B)|_{\m_2}+\ldots+t_r(-B)|_{\m_r}, \ \ t_i>0, \ \ i=1,\ldots,r.$$ 

Let ${X_{\alpha}}$ be a $(-B)$-orthonormal basis adapted to the $\Ad(K)$-invariant decomposition of $\m$, i.e., $X_{\alpha}\in \m_i$ for some $i$, and $\alpha<\beta$ if $i<j$ with $X_{\alpha}\in\m_i$ and $X_{\beta}\in\m_j$. Set $A^{\gamma}_{\alpha\beta}=-B([X_{\alpha},X_{\beta}],X_{\gamma})$, so that $[X_{\alpha},X_{\beta}]_{\m}=\displaystyle\sum_{\gamma}{A^{\gamma}_{\alpha\beta}}X_{\gamma},$ and define 

$$\left[
\begin{array}{cc}
k\\
ij\\
\end{array}
\right]=\sum({A^{\gamma}_{\alpha\beta}})^2,$$ where the sum is taken over all indices $\alpha,\beta,\gamma$ with $X_{\alpha}\in\m_i,X_{\beta}\in\m_j $ and $X_{\gamma}\in\m_k$. Note that $\left[
\begin{array}{cc}
k\\
ij\\
\end{array}
\right]$ is independent of the $(-B)$-orthogonal basis chosen for $\m_i,\m_j,\m_k$, but it depends on the choice of the decomposition of $\m$. In addition, $\left[
\begin{array}{cc}
k\\
ij\\
\end{array}
\right]$ is continuous function on the space of all $(-B)$-orthogonal ordered decomposition of $\m$ into $\Ad(K)-$irreducible summands and also is non-negative and symmetric in all 3 indices. The set $\{X_{\alpha}/\sqrt{t_i};X_{\alpha}\in\m_i\}$ is a orthonormal basis of $\m$ with respect to $$\left\langle \cdot,\cdot\right\rangle=t_1(-B)|_{\m_1}+t_2(-B)|_{\m_2}+\ldots+t_r(-B)|_{\m_r}.$$Then the scalar curvature of $g$ determined by $\left\langle \cdot,\cdot\right\rangle$ is
\begin{equation}
\textrm{scal}(g)=\frac{1}{2}\displaystyle\sum_{l=1}^{r}{\frac{d_l}{t_l}-\frac{1}{4}\displaystyle\sum_{i,j,k}\left[
\begin{array}{cc}
k\\
ij\\
\end{array}
\right]\frac{t_k}{t_it_j}}, \label{curvesc}
\end{equation}
for all $x\in G/K$, $d_i=\dim\m_i$, $i=1,\ldots,r$. See \cite{WZ} for details.

From the above mentioned formula, Nikoronov obtained in his work \cite{nikoronov} formulae for $\text{scal}(g)$ where $g$ is the normal metric on the special cases $W_{1,0}^7$ and $W_{1,1}^7$. We will obtain such a formula for the normal metric on the general case $W_{k,l}^7$ and then, by using the expressions of $\cvg,\cvh,\cvk $ defined on $W_{k,l}^7$, $W_{1,0}^7$ and $W_{1,1}^7$, respectively, we get the following formulae for $\text{scal}(\cvg)$, $\text{scal}(\cvh)$ and $\text{scal}(\cvk)$.


\begin{lemma}
The scalar curvature of the Riemannian homogeneous space  $(W_{k,l}^7=SU(3)/S_{k,l}^1, g_{t_1,t_2, t_3, t_4})$, $\gcd(k,l)=1, k>l>0$, $g_{t_1,t_2, t_3, t_4}$ induced by the inner product $$\left\langle \cdot,\cdot\right\rangle=g_{1\cdot K}=t_1(-B)|_{\m_1}+t_2(-B)|_{\m_2}+ +t_3(-B)|_{\m_3}+t_4(-B)|_{\m_4}$$on $\m$, with $t_i$ positive constants, $B$ (the extension of) the Cartan-killing form of $\mathfrak{su}(3)$  and \linebreak $\m=\m_1\oplus\m_2\oplus\m_3\oplus\m_4$ the isotropy representation of $K=S^1_{k,l}$, is given by  

\begin{eqnarray*}
\emph{scal}(g_{t_1,t_2, t_3, t_4})&=&\dfrac{-(t_1^2+t_2^2+t_3^2)+6(t_1t_2+t_2t_3+t_1t_3)}{6t_1t_2t_3}-\\
                                  & &-\dfrac{t_4}{8\gamma}\left(\dfrac{(k+l)^2}{t_1^2}+\dfrac{l^2}{t_2^2}+\dfrac{k^2}{t_3^2}\right)\\
\end{eqnarray*}

\end{lemma}

 {\bf Proof:} In the case of $W_{k,l}^7=SU(3)/S^1_{k,l}$, $\gcd(k,l)=1, k>l>0$, since the isotropic summands of the decomposition $\m=\m_1\oplus\m_2\oplus\m_3\oplus\m_4,$ are mutually inequivalents, it follows that a $SU(3)$-invariant metric $g$ on $W_{k,l}^7$ is represented by a inner product $$\left\langle \cdot,\cdot\right\rangle=g_{1\cdot K}=t_1(-B)|_{\m_1}+t_2(-B)|_{\m_2}+ +t_3(-B)|_{\m_3}+t_4(-B)|_{\m_4}$$on $\m$, with $t_i$ positive constants, $B$ is (the extension of) the Cartan-killing form of $\mathfrak{su}(3)$  and \linebreak $\m=\m_1\oplus\m_2\oplus\m_3\oplus\m_4$ is the isotropy representation of $K=S^1_{k,l}$, accordingly the construction in Section \ref{isotropia}. The set $\{Y_1=X_1/\sqrt{12t_1,},Y_2=X_2/\sqrt{12t_1},Y_3=X_3/\sqrt{12t_2}, Y_4=X_4/\sqrt{12t_2}, Y_5=X_5/\sqrt{12t_3},Y_6=X_6/\sqrt{12t_3}, Y_7=X_0/\sqrt{12t_4}\}$ is a orthonormal basis of $\m$ with respect to $\left\langle \cdot,\cdot\right\rangle$. 

    It is well known that the expression for the value at some $p_0\in W^7_{k,l}$ of the Ricci tensor field $Ric$ on the unimodular Riemannian homogeneous space $W_{k,l}^7=SU(3)/S^1_{k,l}$ endowed with a $SU(3)$-invariant metric induced by $\left\langle \cdot,\cdot\right\rangle$ is given by $$Ric(Y,Y)=-\dfrac{1}{2}\sum_j\langle[Y,Y_j]_{\m},[Y,Y_j]_{\m} \rangle - \dfrac{1}{2}B(Y,Y)+\dfrac{1}{4}\sum_{i,j}\langle[Y_i,Y_j]_{\m},Y \rangle^2,$$ for $Y\in \m$ and $B$ the Cartan-Killing form of $\mathfrak{su}(3)$ (see \cite{besse}). We have that $B(X,Y)=-6\text{ReTr}(XY)$ and then 

    \begin{eqnarray*}
     Ric(Y_1,Y_1)&=&-\dfrac{1}{2}\left(\dfrac{(k+l)^2t_4}{4\gamma t_1^2}+\dfrac{t_3}{12t_1t_2}+\dfrac{t_3}{12t_1t_2}+\dfrac{t_2}{12t_1t_3}+\dfrac{t_2}{12t_1t_3}+\dfrac{(k+l)^2)}{4\gamma t_4}\right)+\\
                 & &+\dfrac{1}{2t_1}(-\dfrac{1}{2}\text{ReTr}(X_1X_1))+\dfrac{1}{2}\left(\dfrac{(k+l)^2}{4\gamma t_4}+\dfrac{t_1}{12t_2t_3}+\dfrac{t_1}{12t_2t_3}\right) \\
                 &=& \dfrac{t_1^2-t_2^2-t_3^2+6t_2t_3}{12t_1t_2t_3}-\dfrac{t_4(k+l)^2}{8\gamma t_1^2} \\
                 &=& Ric(Y_2,Y_2);
    \end{eqnarray*}
analogously,

    $$ Ric(Y_3,Y_3)= \dfrac{t_2^2-t_3^2-t_1^2+6t_1t_3}{12t_1t_2t_3}-\dfrac{t_4l^2}{8\gamma t_2^2}= Ric(Y_4,Y_4); $$

    $$ Ric(Y_5,Y_5)= \dfrac{t_3^2-t_2^2-t_1^2+6t_1t_2}{12t_1t_2t_3}-\dfrac{t_4k^2}{8\gamma t_3^2}= Ric(Y_6,Y_6); $$
 and finally,                 
    
\begin{eqnarray*}
     Ric(Y_7,Y_7)&=&-\dfrac{1}{2}\left(\dfrac{(k+l)^2}{4\gamma t_4}+\dfrac{(k+l)^2}{4\gamma t_4}+\dfrac{l^2}{4\gamma t_4}+\dfrac{l^2}{4\gamma t_4}+\dfrac{k^2}{4\gamma t_4}+\dfrac{k^2}{4\gamma t_4}\right)+\\
                 & &+\dfrac{1}{2t_4}+\dfrac{1}{2}\left(\dfrac{(k+l)^2t_4}{4\gamma t_1^2}+\dfrac{l^2t_4}{4\gamma t_2^2}+\dfrac{k^2t_4}{4\gamma t_3^2}\right) \\
                 &=&\dfrac{t_4}{8\gamma}\left(\dfrac{(k+l)^2}{t_1^2}+\dfrac{l^2}{t_2^2}+\dfrac{k^2}{t_3^2}\right).\\
    \end{eqnarray*}

The scalar curvature of the $SU(3)$-invariant metric determined by $\left\langle\cdot , \cdot\right\rangle$ is the trace of $Ric$, i.e., $\sum_jRic(Y_j,Y_j)$. Hence, we have the following formula for scalar curvature of such a metric in terms of $t_1,t_2, t_3, t_4$,  

\begin{eqnarray*}
\text{scal}(g_{t_1,t_2, t_3, t_4})&=&\dfrac{-(t_1^2+t_2^2+t_3^2)+6(t_1t_2+t_2t_3+t_1t_3)}{6t_1t_2t_3} -\\
                                  & &-\dfrac{t_4}{8\gamma}\left(\dfrac{(k+l)^2}{t_1^2}+\dfrac{l^2}{t_2^2}+\dfrac{k^2}{t_3^2}\right)\\ \label{scalG}
\end{eqnarray*}
\cqd

\begin{proposition} \emph{[$(\bf {W_{k,l}^7,\cvg)}$]} Let $(W_{k,l}^7,\cvg)$ be the canonical variation of $(W_{k,l}^7,g)$, where $g$ is the normal metric on $W_{k,l}^7$, $\gcd(k,l)=1, k>l>0$. Then, the function $\emph{scal}(t)$, for each $t>0$, 
\begin{equation}
\emph{scal}(t)=-\dfrac{(10\gamma-3k^2)t^2}{24\gamma}+\dfrac{(16\gamma-3k^2)}{24\gamma t^2}+2\label{scal}
\end{equation}
gives the scalar curvature of $\cvg$, $\gamma=k^2+kl+l^2$. 
\end{proposition}
{\bf Proof:} The canonical variation $\cvg$ is given by $${(\cvg)}_{eK}=(-B)|_{\m_1}+(-B)|_{\m_2}+t^2(-B)|_{\m_3}+t^2(-B)|_{\m_4},$$ so that, from the formula of $\text{scal}(g_{t_1,t_2, t_3, t_4})$ with $t_1=t_2=1$ and $t_3=t_4=t^2$ follows the above expression given in \ref{scal}.

\cqd



\begin{proposition} \emph{[$\bf {(W_{1,0}^7,\cvh)}$]} Let $(W_{1,0}^7,\cvh)$ be the canonical variation of $(W_{1,0}^7,g)$, where $g$ is the normal metric on $W_{1,0}^7$. Then, the function $\emph{scal}(t)$, for each $t>0$, 
\begin{equation}
\emph{scal}(t)=\dfrac{-5t^4+24t^2+8}{t^2}\label{scal1}
\end{equation}
gives the scalar curvature of $\cvh$. 
\end{proposition}
{\bf Proof:} In this case, one has 
$${(\cvh)}_{1\cdot SO(2)}=(-B)|_{\mathfrak{p}_1}+(-B)|_{\mathfrak{p}_2}+t^2(-B)|_{\mathfrak{p}_3}+t^2(-B)|_{\mathfrak{p}_4}.$$

We have also that invariant metrics on $W^7_{1,0}=SU(3)/SO(2)$ are given by $$g_{1\cdot SO(2)}=t_1(-B)|_{\mathfrak{p}_1}+t_2(-B)|_{\mathfrak{p}_2}+t_3(-B)|_{\mathfrak{p}_3}+t_4(-B)|_{\mathfrak{p}_4},$$where $\mathfrak{p}_1, \mathfrak{p}_2, \mathfrak{p}_3, \mathfrak{p}_4$ are isotropic summands of $\m$ and $B$ is the Cartan-Killing form of $\mathfrak{su}(3)$, according \ref{10}. Nikoronov in his work \cite{nikoronov} determined the following formula for $\text{scal}(g_{t_1,t_2,t_3,t_4})$ of such a metric on $W^7_{1,0}$ in terms of $t_1,t_2,t_3,t_4$:

\begin{eqnarray*}
\text{scal}(g_{t_1,t_2, t_3, t_4})&=&\dfrac{12}{t_1}+\dfrac{12}{t_2}+\dfrac{12}{t_3}+\dfrac{6a}{t_4}-\dfrac{3-3a}{2}\left(\dfrac{t_4}{t_1^2}+\dfrac{t_4}{t_2^2}\right)-\\
                                  & &-2\left(\dfrac{t_1}{t_2t_3}+\dfrac{t_2}{t_1t_3}+\dfrac{t_3}{t_2t_1}\right)-3a\left(\dfrac{t_1}{t_2t_4}+\dfrac{t_2}{t_1t_4}+\dfrac{t_4}{t_2t_1}\right)\\ 
\end{eqnarray*}
with $a=\sen^2(\alpha)$, $\alpha$ given in \ref{10}. Taking $t_1=t_2=1$ and $t_3=t_4=t^2$ in $\text{scal}(g_{t_1,t_2, t_3, t_4})$, we obtain \ref{scal1}.

\cqd

\begin{proposition} \emph{[$\bf {(W_{1,1}^7,\cvk)}$]} Let $(W_{1,1}^7,\cvk)$ be the canonical variation of $(W_{1,1}^7,g)$, where $g$ is the normal metric on $W_{1,1}^7$. Then, the function $\emph{scal}(t)$, for each $t>0$, 
\begin{equation}
\emph{scal}(t)=-\dfrac{3t^4-24t^2-6}{t^2}\label{scal2}
\end{equation}
gives the scalar curvature of $\cvk$. 
\end{proposition}
{\bf Proof:} The canonical variation $\cvk$ on $W^7_{1,1}=SU(3)/S^1_{1,1}$ is given by
$${(\cvk)}_{eK}=(-B)|_{\mathfrak{p}_1}+(-B)|_{\mathfrak{p}_2}+t^2(-B)|_{\mathfrak{p}_3}+t^2(-B)|_{\mathfrak{p}_4}+t^2(-B)|_{\mathfrak{p}_5}.$$

We have also that invariant metrics on $W^7_{1,1}$ are given by $$g_{1\cdot K}=t_1(-B)|_{\mathfrak{p}_1}+t_2(-B)|_{\mathfrak{p}_2}+t_3(-B)|_{\mathfrak{p}_3}+t_4(-B)|_{\mathfrak{p}_4}+t_5(-B)|_{\mathfrak{p}_5},$$where $\mathfrak{p}_1, \mathfrak{p}_2, \mathfrak{p}_3, \mathfrak{p}_4, \mathfrak{p}_5$ are isotropic summands of $\m$ and $B$ is the Cartan-Killing form of $\mathfrak{su}(3)$, according \ref{10}. In \cite{nikoronov} we also have following formula for $\text{scal}(g_{t_1,t_2,t_3,t_4, t_5})$ of such a metric on $W^7_{1,1}$ in terms of $t_1,t_2,t_3,t_4, t_5$:

\begin{eqnarray*}
\text{scal}(g_{t_1,t_2, t_3, t_4,t_5})&=&\dfrac{12}{t_1}+\dfrac{12}{t_2}+\dfrac{6}{t_3}+\dfrac{6}{t_4}+\dfrac{6}{t_5}-2\left(\dfrac{t_3}{t_4t_5}+\dfrac{t_4}{t_3t_5}+\dfrac{t_5}{t_3t_4}\right)-\\
                                  & &-\dfrac{1-b^2}{2}\left(\dfrac{4}{t_3}+\dfrac{t_3}{t_1^2}+\dfrac{t_3}{t_2^2}\right)-\dfrac{b^2-a^2}{2}\left(\dfrac{4}{t_4}+\dfrac{t_4}{t_1^2}+\dfrac{t_2}{t_3t_1}\right)-\\ 
                                  & &-\dfrac{a^2}{2}\left(\dfrac{4}{t_5}+\dfrac{t_5}{t_1^2}+\dfrac{t_5}{t_2^2}\right)-b^2\left(\dfrac{t_3}{t_1t
                                  _2}+\dfrac{t_1}{t_2t
                                  _3}+\dfrac{t_2}{t_3t_1}\right)-\\
                                  & &-(1+a^2-b^2)\left(\dfrac{t_4}{t_1t_2}+\dfrac{t_1}{t_2t
                                  _4}+\dfrac{t_2}{t_4t_1}\right)-\\
                                  & &-(1-a^2)\left(\dfrac{t_3}{t_1t_2}+\dfrac{t_1}{t_2t
                                  _5}+\dfrac{t_2}{t_5t_1}\right)\\
\end{eqnarray*}
with $a=\cos\varphi$, $b=\cos\varphi/Q$, $\varphi$, $Q$ given in \ref{11}. Taking $t_1=t_2=1$ and $t_3=t_4=t_5=t^2$ in $\text{scal}(g_{t_1,t_2, t_3, t_4,t_5})$, we obtain \ref{scal2}.

\cqd

%

\subsection{Bifurcation and Local Rigidity}

If the Morse index changes when passing a degeneracy instant this is actually a bifurcation instant. We will apply such principle to find bifurcation instants for the canonical variations $\cvg,\cvh,\cvk,\cvm,\cvn, 0<t<1$. 

It was defined previously that a {\it degeneracy instant} $t_{\ast}>0$ for $g_t$ in ${\mathcal{R}}^k(M)$, with $g_1=g$, is an instant such that $\dfrac{\textrm{scal}(g_{t_{\ast}})}{m-1}\in\sigma(\Delta_{t_{\ast}})$.

Denoting by $$\{0<\lambda^g_1\leq \lambda^g_2\leq\ldots\leq\lambda^g_j\leq\ldots\}$$ the sequence of positive eigenvalues of $\Delta_g$, the Morse index of a Riemannian metric $g$ is $$N(g)=\max\left\{j\in\mathbb{N};\lambda^g_j<\dfrac{\textrm{scal}(g)}{m-1}\right\},$$ where is the Laplacian $\Delta_g$ acting on $C^{\infty}(M)$, $M$ provided with the Riemannian metric $g$ and $m=\dim M.$ 

The following result is a sufficient condition for a degeneracy instant $0<t_{\ast}<1$ to be a bifurcation instant when $\dfrac{\text{scal}(t_{\ast})}{m-1}$ is a constant eigenvalue of the Laplacian $\Delta_{t_{\ast}}$. 

\begin{proposition}[\cite{pacificjournal}]\label{propbifurcation} Let $(M,g)$ be a closed Riemannian manifold with $\dim M\geq 3$ and \linebreak $\pi:(M,g)\longrightarrow (B,h)$ a Riemannian submersion with totally geodesic fibers isometric to $(F,\kappa)$, where $\dim F\geq 2$ and $\emph{scal}(F)>0$. Denote by $\lambda\in\sigma(\Delta_h)\subset\sigma(\Delta_g)$ a constant eigenvalue of $\Delta_{t_{\ast}}$ such that $\dfrac{\emph{scal}(g_{t_{\ast}})}{m-1}=\lambda$, $g_{t_{\ast}}$ canonical variation of $g$ at $0<t_{\ast}<1$ and $\Delta_{t_{\ast}}$ the Laplacian on $(M,g_{t_{\ast}})$. If $$\dfrac{\emph{scal}(g_{t_{\ast}})}{m-1}<\lambda^{1,1}(t_{\ast})=\mu_1+(\frac{1}{t_{\ast}^2}-1)\phi_1,$$ $\mu_1\in \sigma(\Delta_g)$ the first positive eigenvalue of $\Delta_g$ and $\phi_1\in \sigma(\Delta_v)$ the first positive eigenvalue of the vertical Laplacian $\Delta_v$, then $t_{\ast}$ is a bifurcation instant for $g_t$.
\end{proposition}

\begin{lemma}\label{autoconstante}Every degeneracy instant $t_{\ast}$ for $\cvg,\cvh,\cvk$ is such that $$\dfrac{\emph{scal}(t_{\ast})}{m-1}\in\sigma(\Delta_{\breve{g}})\subset \sigma(\Delta_{t_{\ast}}),$$ with $\sigma(\Delta_{t_{\ast}})$ being the spectrum of the Laplacians on the respective total spaces $W^7_{k,l}, W^7_{1,0}, W^7_{1,1} $ and $\sigma(\Delta_{\breve{g}})$ being the spectrum on the correspondent basis spaces $(G/H,\breve{g})$, $m=\dim W^7_{k,l}=7.$ In other words, $\dfrac{\emph{scal}({t_{\ast}})}{m-1}$ is eigenvalue of $\Delta_{t_{\ast}}$ if and only if $\dfrac{\emph{scal}({t_{\ast}})}{m-1}$ is a constant eigenvalue $\lambda^{r,0}(t)\in\sigma(\Delta_{t_{\ast}})$, for some $1\leq r\in\mathbb{Z}.$

\end{lemma}

{\bf Proof:} It is known that an eigenvalue of $\Delta_t$ can be written as $$\lambda^{r,j}(t)=\mu_r+\left(\frac{1}{t^2}-1\right)\phi_j,$$ for some $\mu_r\in\sigma(\Delta_g)$ and $\phi_j\in\sigma(\Delta_{\hat{g}})$. To proof that $\dfrac{\text{scal}({t_{\ast}})}{m-1}$ is a constant eigenvalue $\lambda^{r,0}(t)\in\sigma(\Delta_{t_{\ast}})$, for some $1\leq r\in\mathbb{Z},$ it is enough to show that there are no $r,j>0$ such that $$\dfrac{\text{scal}({t_{\ast}})}{m-1}=\lambda^{r,j}(t)$$ nor $j>0$ such that $$\dfrac{\textrm{scal}(t)}{m-1}=\lambda^{0,j}(t).$$
In order to prove the inequality $$\dfrac{\text{scal}(t)}{m-1}<\lambda^{1,1}(t)=\mu_1+\left(\dfrac{1}{t^2}-1\right)\phi_1,\forall\, 0<t\leq 1,$$ we define for each canonical variation $\cvg,\cvh,\cvk$ the function $f$ given by $$f(t)=\dfrac{\text{scal}(t)}{6}-\lambda^{1,1}(t)=\dfrac{\text{scal}(t)}{6}-\mu_1-\left(\dfrac{1}{t^2}-1\right)\phi_1,\forall\, t>0,$$and we verify that $f$ is strictly negative for $0<t\leq 1$. Indeed, in the case of the canonical variation $\cvg$ on $W^7_{k,l}$, $\gcd(k,l)=1$ and $k>l>0$, we obtain $$f(t)=-\dfrac{(10\gamma-3k^2)t^2}{144\gamma}+\dfrac{(16\gamma-3k^2)t^2}{144\gamma t^2}+\dfrac{1}{3}-1-\left(\dfrac{1}{t^2}-1\right)8, $$since $\mu_1=1$ and $\phi=8$. It follows that $f(1)=-\dfrac{5}{8}<0$ and 

$$\dfrac{df}{dt}(t)=\dfrac{1139\gamma-(10\gamma-3k^2)t^4}{72\gamma t^3}>0$$
for $0<t\leq 1$. In the special cases, the proof is analogous.
Furthermore, it can be shown, by applying elementary differential calculus, that there is no $0<t<1$ such that $\dfrac{\textrm{scal}(t)}{6}=\lambda^{0,j}(t),$ for all integer $j\geq 1$. We will present the proof of this property for the canonical variation $\cvg$. The proof of the another cases are entirely analogous.

From the description of the spectra of Aloff-Wallach and the complex projective space $\mathbb{CP}^2$, the first positive eigenvalue of the Laplacian $\Delta_t$ acting on functions on the total space $(W^7_{k,l},\cvg)$ is $\lambda_1(t)=1$, for all $0<t\leq 1$. Since $$\dfrac{\textrm{scal}(t)}{m-1}=\dfrac{\textrm{scal}(t)}{6}=-\dfrac{(10\gamma-3k^2)t^2}{144\gamma}+\dfrac{(16\gamma-3k^2)}{144\gamma t^2}+\dfrac{1}{3}$$ for all $t>0$, we have 
\begin{eqnarray*}
\dfrac{\textrm{scal}(g_t)}{m-1}&=&-\dfrac{(10\gamma-3k^2)t^2}{144\gamma}+\dfrac{(16\gamma-3k^2)}{144\gamma t^2}+\dfrac{1}{3}\\                            & < &\lambda_1(t)=1\\
&\Leftrightarrow & \sqrt{\dfrac{-96\gamma+\sqrt{(96\gamma)^2+4(16\gamma-3k^2)(10\gamma-3k^2)}}{2(10\gamma-3k^2)}}<t\leq 1, 
\end{eqnarray*}
for all integers $k,l$ such that $0<l<k$.
Thus, we have that $\cvg$ is locally rigidity at all instant 
$$t_{\ast}\in \ \  ]\sqrt{\dfrac{-96\gamma+\sqrt{(96\gamma)^2+4(16\gamma-3k^2)(10\gamma-3k^2)}}{2(10\gamma-3k^2)}},1].$$ Hence, there are no degeneracy instants for $\cvg$ in the interval $]b,1],$ with $$b=\sqrt{\dfrac{-96\gamma+\sqrt{(96\gamma)^2+4(16\gamma-3k^2)(10\gamma-3k^2)}}{2(10\gamma-3k^2)}}.$$

The fiber of the canonical fibration $(W^7_{k,l},\cvg)$ is the homogeneous lens space $(U(2)/S^1_{k,l},\hat{g})$, with the induced metric $\hat{g}$, represented by the inner product $(-B)|_{\p}$. The first positive eigenvalue $\phi_1$ of the Laplacian $\Delta_{\hat{g}}$ of such homogeneous metric is equal to 8, i.e., if $\phi\in\sigma(\Delta_{\hat{g}})$, then $\phi\geq 8$.

Define $\varphi_c:]0,1[\longrightarrow \R$ by $$\varphi_c(t)=\dfrac{\textrm{scal}(g_t)}{m-1}-\left(\dfrac{1}{t^2}-1\right)\cdot c,$$ for some fixed $c\geq 8$. We have that $$\dfrac{d}{dt}(\varphi_c(t))=-\frac{16\gamma-3k^2}{72t^3\gamma}+\frac{t\left(3k^2-10\gamma\right)}{72\gamma}+\frac{2c}{t^3}>0,$$and $$\varphi_c(b)=\frac{3k^2c+(c+1)\sqrt{2395k^4+4850k^3l+7314k^2l^2+4928kl^3+2464l^4}-2(29c+24)\gamma}{\sqrt{2395k^4+4850k^3l+7314k^2l^2+4928kl^3+2464l^4}-48\gamma}<0,$$$\forall \, 0<t<1, c\geq 8,$ and for any integers $0<k<l$. It follows that $\varphi_c(t)$ is negative for all $0<t<b$ and $$\dfrac{\textrm{scal}(g_t)}{m-1}<\left(\dfrac{1}{t^2}-1\right)\cdot c,$$ for any $c\geq 8$, that is, there is no $0<t<1$ such that $\dfrac{\textrm{scal}(t)}{m-1}=\lambda^{0,j}(t)=\left(\dfrac{1}{t^2}-1\right)\cdot\phi_j,$ for all $\phi_j\in\sigma(\Delta_{\hat{g}})$ and we complete the proof for the case of the canonical variation $\cvg$.

As well as for the case of $\cvg$, define for $\cvh, \cvk,$ the function $\varphi_c(t)$ as above with $c\geq \phi_1$, where $\phi_1$ is the first positive eigenvalues of the Laplacians on their respective fibers. The statement follows by checking that $\dfrac{d}{dt}(\varphi_c(t))>0, \forall\, 0<t<1$ and $\varphi_c(b)<0$, where $0<b<1$ is such that $$\dfrac{\textrm{scal}(g_t)}{m-1}<\lambda_1(t), \forall\, b<t<1,$$with $\lambda_1(t)$ first positive eigenvalues of the Laplacians.
 \cqd
 


The degeneracy instants for the conical variation $\cvg$, $0<t<1$, on $W^7_{k,l}$ are given in our following theorem.

\begin{theorem}\label{teo-princ01}
Let $\cvg$ be the aforementioned canonical variation on $W^7_{k,l}$ and take $$b=\sqrt{\dfrac{-96\gamma+\sqrt{(96\gamma)^2+4(16\gamma-3k^2)(10\gamma-3k^2)}}{2(10\gamma-3k^2)}}.$$ Thus, the degeneracy instants for $\cvg$ in $]0,1[$ form a decreasing sequence $\{t^{\text{\bf g}}_{q}\}\subset \, ]0,b]$ such that $t^{\text{\bf g}}_{q}\rightarrow 0$ as $q\rightarrow 0$, with $t^{\text{\bf g}}_1=b$ and for $q>1$, 
\begin{equation}
t^{\text{\bf g}}_{q}=\dfrac{1}{\sqrt{2(10\gamma-3k^2})}\sqrt{48\gamma(1-2q -q^2)+\sqrt{(48\gamma)^2(1-2q -q^2)^2-4(3k^2-10\gamma)(16\gamma -3k^2)}}, \label{seqsu}
\end{equation}
$0<l<k$, $k,l$ relatively prime and $\gamma=k^2+kl+l^2.$
\end{theorem}

{\bf Proof:} We must determine all $t\in \ \ ]0,b]$ such that $\dfrac{\textrm{scal}(t)}{m-1}\in\sigma(\Delta_t)$, since by the previous theorem $g_t$ is locally rigidity for $b<t<1$. From Lemma \ref{autoconstante}, if $\dfrac{\textrm{scal}(t)}{m-1}$ is an eigenvalue of $\Delta_t$, then $\dfrac{\textrm{scal}(t)}{m-1}\in\Delta_{\breve{g}}$.  

Thus, it remains to verify for which instants $0<t<b<1$ one has $\dfrac{\textrm{scal}(t)}{m-1}=\lambda^{r,0}(t)$.  

We have that $\lambda^{r,0}(t)$ is eigenvalue of $\Delta_t$ if and only if $\lambda^{r,0}(t)$ belongs to the spectrum of the Laplacian on the basis $\mathbb{CP}^2$, provided with the symmetric metric $\breve{g}$ represented by the inner product $(\cdot,\cdot)=-B|_{\mathfrak{q}}$, where $\mathfrak{q}$ is the horizontal space of the original fibration. The spectrum of the Laplacian $\Delta_{\breve{g}}$ on $(\mathbb{CP}^n,\breve{g})$ is $$ \sigma(\Delta_{\breve{g}})=\{\beta_q=\dfrac{q(q+2)}{3};q\in\mathbb{N}\}.$$
 
The degeneracy instants for $\cvg$ in $]0,b]$ are the real values $t^g_q$, solutions of the equation $$\dfrac{\textrm{scal}(t)}{m-1}=\beta_q=\dfrac{q(q+2)}{3}. $$ The explicit formula for $t^{\text{\bf g}}_{q}$, which represents the solution of the above equation in $t$, is presented in \eqref{seqsu}. Note that the constant eigenvalues $\beta_q=\dfrac{q(q+2)}{3}$ tend to $+\infty$ as $q\rightarrow \infty$ and, since $\dfrac{\text{scal}(t)}{m-1}$ is continuous and tends to $+\infty$ as $t\rightarrow 0$, $t^{\text{\bf g}}_{q}\rightarrow 0.$ \cqd 


We also compute the Morse index of each $\cvg$ as a critical point of the Hilbert-Einstein functional.

The degeneracy instants for the canonical variation $\cvh$, $0<t\leq 1$, on the special case $W^2_{1,0}$ are given in our following theorem.
\begin{theorem} \label{teo-princ-02}
Let $\cvh$ be the above canonical variation on $W^7_{1,0}$ and take $b=2.$ Thus, the degeneracy instants for $\cvh$ in $]0,b]$ form a decreasing sequence $\{t^{\text{\bf h}}_{q}\}\subset \, ]0,b]$ such that $t^{\text{\bf h}}_{q}\rightarrow 0$ as $q\rightarrow 0$, with $t^{\text{\bf h}}_1=b$ and for $q>1$, 
\begin{equation}
t^{\text{\bf h}}_{q}=\frac{\sqrt{24-4q-2q^2+\sqrt{\left(-2q^2-4q+24\right)^2+160}}}{\sqrt{10}}. \label{seqso1}
\end{equation}

\end{theorem}

{\bf Proof:} We must determine all $t\in \ \ ]0,b]$ such that $\dfrac{\textrm{scal}(t)}{6}\in\sigma(\Delta_t)$, since by the previous theorem $\cvh$ is locally rigidity for $b<t$. From Lemma \ref{autoconstante}, if $\dfrac{\textrm{scal}(t)}{6}$ is an eigenvalue of $\Delta_t$, then $\dfrac{\textrm{scal}(t)}{6}\in\Delta_{\breve{g}}$.

Thus, it remains to verify for which instants $0<t<b<1$ one has $\dfrac{\textrm{scal}(t)}{6}=\lambda^{r,0}(t)$. 

From Proposition \ref{constanteigen}, Section \ref{Section2.2}, $\lambda^{r,0}(t)$ is eigenvalue of $\Delta_t$ if and only if $\lambda^{r,0}(t)$ belongs to the spectrum of the Laplacian on the basis $\mathbb{CP}^2$, provided with the symmetric metric $\breve{g}$ represented by the inner product $(\cdot,\cdot)=-B|_{\mathfrak{q}}$, $\mathfrak{q}$ horizontal space of the original fibration. The spectrum of the Laplacian $\Delta_{\breve{g}}$ on $(\mathbb{CP}^2,\breve{g})$ is $$ \sigma(\Delta_{\breve{g}})=\{\beta_q=\dfrac{q(q+2)}{3};q\in\mathbb{N}\}.$$
 
The degeneracy instants for $\cvh$ in $]0,b]$ are the real values $t^{\text{\bf h}}_{q}$, solutions of the equation $$\dfrac{\textrm{scal}(t)}{6}=\beta_q=\dfrac{q(q+2)}{3}. $$ The explicit formula for $t^{\text{\bf h}}_q$, which represents the solution of the above equation in $t$, is presented in \eqref{seqso1}. Note that the constant eigenvalues $\beta_q=\dfrac{q(q+2)}{3}$ tend to $+\infty$ as $q\rightarrow \infty$ and, since $\dfrac{\text{scal}(t)}{6}$ is continuous and tends to $+\infty$ as $t\rightarrow 0$, $t^{\text{\bf h}}_{q}\rightarrow 0.$ \cqd 

\begin{theorem} \label{teo-princ03}
Let $\cvk$ be the above canonical variation on $W^7_{1,1}$ and take $b=\sqrt{3+\sqrt{10}}.$ Thus, the degeneracy instants for $\cvk$ in $]0,b]$ form a decreasing sequence $\{t^{\text{\bf k}}_{q}\}\subset \, ]0,b]$ such that $t^{\text{\bf k}}_{q}\rightarrow 0$ as $q\rightarrow 0$, with $t^{\text{\bf k}}_1=b$ and for $q>1$, 
\begin{equation}
t^{\text{\bf k}}_{q}=\frac{\sqrt{24-4q-2q^2+\sqrt{\left(-2q^2-4q+24\right)^2+36}}}{\sqrt{6}}. \label{seqso1}
\end{equation}

\end{theorem}

{\bf Proof:} We must determine all $t\in \ \ ]0,b]$ such that $\dfrac{\textrm{scal}(t)}{6}\in\sigma(\Delta_t)$, since by the previous theorem $\cvk$ is locally rigidity for $b<t$. From Lemma \ref{autoconstante}, if $\dfrac{\textrm{scal}(t)}{6}$ is an eigenvalue of $\Delta_t$, then $\dfrac{\textrm{scal}(t)}{6}\in\Delta_{\breve{g}}$.

Thus, it remains to verify for which instants $0<t<b<1$ one has $\dfrac{\textrm{scal}(t)}{6}=\lambda^{r,0}(t)$. 

From Proposition \ref{constanteigen}, Section \ref{Section2.2}, $\lambda^{r,0}(t)$ is eigenvalue of $\Delta_t$ if and only if $\lambda^{r,0}(t)$ belongs to the spectrum of the Laplacian on the basis $\mathbb{CP}^2$, provided with the symmetric metric $\breve{g}$ represented by the inner product $(\cdot,\cdot)=-B|_{\mathfrak{q}}$, $\mathfrak{q}$ horizontal space of the original fibration. The spectrum of the Laplacian $\Delta_{\breve{g}}$ on $(\mathbb{CP}^2,\breve{g})$ is $$ \sigma(\Delta_{\breve{g}})=\{\beta_q=\dfrac{q(q+2)}{3};q\in\mathbb{N}\}.$$
 
The degeneracy instants for $\cvk$ in $]0,b]$ are the real values $t^{\text{\bf k}}_{q}$, solutions of the equation $$\dfrac{\textrm{scal}(t)}{6}=\beta_q=\dfrac{q(q+2)}{3}. $$ The explicit formula for $t^{\text{\bf k}}_q$, which represents the solution of the above equation in $t$, is presented in \eqref{seqso1}. Note that the constant eigenvalues $\beta_q=\dfrac{q(q+2)}{3}$ tend to $+\infty$ as $q\rightarrow \infty$ and, since $\dfrac{\text{scal}(t)}{6}$ is continuous and tends to $+\infty$ as $t\rightarrow 0$, $t^{\text{\bf k}}_{q}\rightarrow 0.$ \cqd 

We can determine the Morse index $N(t)$ for each $\cvg$, $\cvh$ and $\cvk$ as following.

\begin{proposition}\label{morseindexsu}
The Morse index of $\cvg$ is given by 
\begin{equation} 
N(t)=\left\{
\begin{array}{rcl}
\displaystyle\Sigma_{q=1}^r m_q(\mathbb{CP}^2),& \mbox{if}&t^{\text{\bf g}}_{r+1}\leq t < t^{\text{\bf g}}_{r}\\ 
0, &\mbox{if}& t^{\text{\bf g}}_1\leq t\leq 1
\end{array}
\right.,
\end{equation}where $m_q(\mathbb{CP}^2)$ is the multiplicity of the $q$th eigenvalue of the basis $\mathbb{CP}^2.$ 
\end{proposition}

{\bf Proof:} We established in Lemma \ref{autoconstante} that $\dfrac{\textrm{scal}(t)}{6}\leq\lambda_1(t)$ for all $t\in [b,1]$, with $b=t^{\text{\bf g}}_1,$ so that there are no eigenvalues of $\Delta_t$ that are less than $\dfrac{\textrm{scal}(t)}{6}$ in the interval $[b,1]$. Hence, $N(t)=0$ for $t\in [b,1]$. When $t\rightarrow 0$, whenever $t$ crosses a degeneracy instant $t^{\text{\bf g}}_{q}$, the constant eigenvalue $\lambda^{q,0}(t)$ becomes smaller than $\dfrac{\textrm{scal}(t)}{6}$. 

Therefore, the Morse index increases by the multiplicity of $\lambda^{q,0}(t),$ which is the dimension of the corresponding  eigenspace $E^{0}_{q}$. This dimension is given by the Weyl dimension formula and is positive. Thus, $\dim E^{0}_{q}$ also coincides with the multiplicity of the complex projective space $\mathbb{CP}^2,$ concluding the proof. \cqd

\begin{proposition}\label{morseindexso}
The Morse index of $\cvh$ is given by 
\begin{equation} 
N(t)=\left\{
\begin{array}{rcl}
\displaystyle\Sigma_{q=1}^r m_q(\mathbb{CP}^2),& \mbox{if}&t^{\text{\bf h}}_{r+1}\leq t < t^{\text{\bf h}}_{r}\\ 
0, &\mbox{if}& t^{\text{\bf h}}_1\leq t\leq 1
\end{array}
\right.,
\end{equation}where $m_q(\mathbb{CP}^2)$ is the multiplicity of the $q$th eigenvalue of the basis $\mathbb{CP}^2.$ 
\end{proposition}

\begin{proposition}\label{morseindexso}
The Morse index of $\cvk$ is given by 
\begin{equation} 
N(t)=\left\{
\begin{array}{rcl}
\displaystyle\Sigma_{q=1}^r m_q(\mathbb{CP}^2),& \mbox{if}&t^{\text{\bf k}}_{r+1}\leq t < t^{\text{\bf k}}_{r}\\ 
0, &\mbox{if}& t^{\text{\bf k}}_1\leq t\leq 1
\end{array}
\right.,
\end{equation}where $m_q(\mathbb{CP}^2)$ is the multiplicity of the $q$th eigenvalue of the basis $\mathbb{CP}^2.$ 
\end{proposition}

As consequence of the above results, we have obtained the bifurcation and rigidity instants for $(W^2_{k,l},\cvg)$, $(W^2_{1,0},\cvh), (W^2_{1,1},\cvk)$, $0<t\leq b,$ in the following theorem.

\begin{theorem} For the canonical variations $(W^2_{k,l},\cvg)$, $(W^2_{1,0},\cvh), (W^2_{1,1},\cvk)$, constructed above, the elements of the sequences $\{t^{\text{\bf g}}_{q}\},\{t^{\text{\bf h}}_{q}\}\subset \, ]0,b]$, of degeneracy instants for $\cvg$, $\cvh$ and $\cvk$, given above, are bifurcation instants for  $\cvg$, $\cvh$ and $\cvk$, respectively. Moreover,  $\cvg$, $\cvh$ and $\cvk$ are locally rigid for all $t\in]0,b]\: \setminus \: \{t^{\text{\bf g}}_{q}\}$(respectively $ t\in]0,b]\: \setminus \: \{t^{\text{\bf h}}_{q}\} \, and \, t\in]0,b]\: \setminus \: \{t^{\text{\bf k}}_{q}\})$.
\end{theorem}

{\bf Proof:} It is known that, if $t_{\ast}>0$ is not a degeneracy instant for $\cvg$, $\cvh$ and $\cvk$, then $\cvg$, $\cvh$ and $\cvk$ are locally rigidity at $t_{\ast}.$ In the previous theorems, we proved that the degeneracy instants for $\cvg$, $\cvh$ and $\cvk$ form the sequences $t^{\text{\bf g}}_{q},t^{\text{\bf h}}_{q}$, $t^{\text{\bf k}}_{q}$ and thus, if $t\notin \{t^{\text{\bf g}}_{q}\}$,  $t\notin \{t^{\text{\bf h}}_{q}\}$ or $t\notin \{t^{\text{\bf k}}_{q}\}$ , $t$ must be a local rigidity instant for $\cvg$, $\cvh$ and $\cvk$, respectively. 

The fact that each $t^{\text{\bf g}}_q$, $t^{\text{\bf h}}_q$ or $t^{\text{\bf k}}_q$ are a bifurcation instants follows from Proposition \ref{propbifurcation} and Lemma \ref{autoconstante}. \cqd



\section{Multiplicity of Solutions to the Yamabe Problem}\label{Section3.3}

We now explain how to obtain multiplicity results for the canonical variations $\cvg,\cvh,\cvk$ applying the next proposition due Bettiol and Piccione in \cite{pacificjournal}. We have been interested in determine which conformal classes carry multiple unit volume metrics with constant scalar curvature. 

\begin{proposition}[\cite{pacificjournal}] Let $g_t$, with $t\in \, ]0,b[$, be a family  of metrics on $M$ with $N(t)=N(g_t)>0$, $N(t)$ the Morse index of $g_t,$ and suppose there exists a sequence $\{t_q\}$ in $]0,b[$, that converges to 0, of bifurcation values for $g_t$. Then, there is an infinite subset $\mathcal{G}\subset ]0,b[$ accumulating at $0$, such that for each $t\in\mathcal{G}$, there are at least $3$ solutions to the Yamabe problem in the conformal class $[g_t].$
\end{proposition}

Applying the last Proposition and our bifurcation results, we can determine a lower bound for the number of unit volume metrics with constant scalar curvature in each conformal class $[\cvg], [\cvh], [\cvk]$, respectively, for instants in a given subset $\mathcal{G}\subset ]0,1[$. 

\begin{theorem}
Let $\cvg,\cvh,\cvk$ be the families of homogeneous metrics described above. Then, there exists, for each of such family, a subset $\mathcal{G}\subset ]0,1[$, accumulating at $0$, such that for each $t\in\mathcal{G}$, there are at least $3$ solutions to the Yamabe problem in each conformal class $[\cvg], [\cvh], [\cvk].$
\end{theorem}

{\bf Proof:} It is only necessary to verify that there exists $0<b<1$ such that $N(t)>0$ in $]0,b[$. For each canonical variation $\cvg,\cvh,\cvk$, take a positive number $0<b<1$ such that $\dfrac{\text{scal(t)}}{m-1}<\lambda_1(t)$ for all $t\in\, ]b,1[$, with $m$ denoting the dimension of the respective total space. This implies that $N(t)=0$, for all $t\in\, ]b,1[$, since there are no eigenvalues of $\Delta_t$ less than  
$\dfrac{\text{scal(t)}}{m-1}$. If $t=b$, one has $\dfrac{\text{scal(t)}}{m-1}=\beta_1,$ where $\beta_1$ is the first positive eigenvalue of the Laplacian on the basis. We proved that $t=b$ is a bifurcation instant and then the Morse index changes from $0$ to a positive integer. Hence, for $t\in\, ]0,b[\subset]0,1[,$ we have $N(t)\geq N(b-\epsilon)>0$, since by definition, the Morse index is equal to the number (counting multiplicity) of positive eigenvalues that are less than $\dfrac{\text{scal}(t)}{m-1}$, which is strictly decreasing for $0<t<1$ and $\dfrac{\text{scal}(t)}{m-1}\rightarrow \infty $ as $t\rightarrow 0.$, for $\cvg,\cvh,\cvk$. \cqd

{\bf Acknowledgment:} L. Grama is partially supported by the grants 2018/13481-0 and 2023/13131-8 (FAPESP). K. Lima is supported by UFCG-CCT(2023- 2024).

\bibliography
{unsrt}  


\end{document}